\documentclass[number,citesort,MSNbibl,dvips]{arxbj}
\usepackage[left]{vmkcol}
\usepackage{upgreek}
\usepackage{graphicx}

\aid{0}
\volume{18}
\issue{1}
\pubyear{2012}
\firstpage{322}
\lastpage{359}
\doi{10.3150/10-BEJ334}

\makeatletter
\newcolumntype{d}[1]{D{.}{.}{#1}}
\newtheorem{theorem}{Theorem}
\newtheorem{lemma}{Lemma}
\newtheorem{proposition}{Proposition}
\newremark{remark}{Remark}
\newcommand{\eqref}[1]{(\ref{#1})}
\newcommand{\cal}{\mathcal}

\newcommand{\trace}{\operatorname{tr}}
\newcommand{\Ai}{\mathrm{Ai}}
\renewcommand{\aa}{\alpha}
\newcommand{\zz}{\zeta}
\newcommand{\ww}{\omega}
\newcommand{\kk}{\kappa}
\newcommand{\tee}{\tilde{\varepsilon}}
\newcommand{\tmu}{\tilde{\mu}}
\newcommand{\tsigma}{\tilde{\sigma}}
\newcommand{\MM}{\mathcal{M}}
\newcommand{\EE}{\mathcal{E}}
\newcommand{\NN}{\mathcal{N}}
\renewcommand{\to}{\rightarrow}
\newcommand{\Oh}[1]{{\mathrm{O}({#1})}}
\newcommand{\convd}{\stackrel{\mathcal D}{\longrightarrow}}
\newcommand\id{{\mathbf1}}

\newcommand{\tauexpr}[1]{\mu_{n,N} + {#1}\sigma_{n,N}}
\newcommand{\xnexpr}[1]{\tmu_{n,N} + {#1}\tsigma_{n,N}}

\makeatother

\begin{document}
\begin{frontmatter}

\title{Accuracy of the Tracy--Widom limits for the extreme eigenvalues in white Wishart matrices}
\runtitle{Extreme eigenvalues in white Wishart matrices}

\begin{aug}
\author{\fnms{Zongming} \snm{Ma}\corref{}\ead[label=e1]{zongming@wharton.upenn.edu}}
\runauthor{Z. Ma}
\address{Department of Statistics, The Wharton School, University of
Pennsylvania, Philadelphia, PA 19104, USA.
\printead{e1}}
\end{aug}

\received{\smonth{4} \syear{2009}}
\revised{\smonth{8} \syear{2010}}

%
\begin{abstract}
The distributions of the largest and the smallest eigenvalues of a
$p$-variate sample covariance matrix $S$ are of great importance in
statistics. Focusing on the null case where $nS$ follows the standard
Wishart distribution $W_p(I,n)$, we study the accuracy of their scaling
limits under the setting: $n/p\to\gamma\in(0,\infty)$ as $n\to
\infty$.
The limits here are the orthogonal Tracy--Widom law and its reflection about
the origin.

With carefully chosen rescaling constants, the approximation to the rescaled
largest eigenvalue distribution by the limit attains accuracy of order
$\Oh{\min(n,p)^{-2/3}}$. If $\gamma> 1$, the same order of accuracy is
obtained for the smallest eigenvalue after incorporating an additional log
transform. Numerical results show that the relative error of approximation
at conventional significance levels is reduced by over $50\%$ in rectangular
and over $75\%$ in `thin' data matrix settings, even with
$\min(n,p)$ as small as $2$.
\end{abstract}

%

\end{frontmatter}

\section{Introduction}
\label{sec:introduction}

Understanding the behavior of the extreme eigenvalues of a sample
covariance matrix~$S$ is important in a large number of multivariate
statistical problems. As an example, consider one of the most common
inference problems: testing the null hypothesis that the population
covariance is identity. Roy's union intersection principle
\cite{roy53} suggests that we reject the null hypothesis for large
values of the largest eigenvalue of $S$ (or for small values of the
smallest eigenvalue). Naturally, the next question is: How should the
$p$-value be calculated?

To address this issue, and many others, it is necessary to examine the null
distributions of the extreme sample eigenvalues. In this paper, we restrict
ourselves to the Gaussian framework. In particular, let $X$ be an
$n\times p$
data matrix whose row vectors are i.i.d. samples from the $N_p(0, I)$
distribution. The $p\times p$ matrix $A = X'X$ then follows a standard Wishart
distribution: $A\sim W_p(I,n)$, and is called a (real) white Wishart matrix.
The ordered eigenvalues of $A$ are denoted by $\lambda_1 \geq\cdots
\geq
\lambda_p$. Our interest lies in $\lambda_1$ and $\lambda_p$, as $A
= nS$.

The exact evaluation of the marginal distributions of these eigenvalues
is difficult, even in the null case considered here. See, for example,
Muirhead \cite{muirhead}, Section 9.7. An alternative approach is to
approximate them by their asymptotic limits. For the problem we are
concerned with, Anderson \cite{anderson}, Chapter 13, summarized the
classical results under the conventional asymptotic regime: $p$ holds
fixed and $n$ tends to infinity.

However, for a wide range of modern data sets (microarray data, stock
prices, weather forecasting, etc.), the number of features $p$ is very
large while the number of observations~$n$ is much smaller than or just
comparable to $p$. For these situations, the classical asymptotics is
not always appropriate and different asymptotic theories are
needed. Borrowing tools from random matrix theory, especially those
established by Tracy and Widom \cite{tw94,tw96,tw98}, Johnstone
\cite{johnstone01} showed that under the asymptotic regime
%
%
\begin{equation}
\label{eq:asymptotics}
p \to\infty, \qquad  n = n(p)\to\infty \quad \mbox{and} \quad  n/p\to\gamma\in
(0,\infty),
\end{equation}
the largest eigenvalue $\lambda_1$ in $A$ has the weak limit
%
%
\begin{equation}
\label{eq:loe-lim}
\frac{\lambda_1 - \mu_{p}}{\sigma_{p}} \convd F_1 ,
\end{equation}
where the centering and scaling constants are defined as
%
%
\begin{equation}
\label{eq:orig-center-scale}
\mu_{p} = \bigl(\sqrt{n-1} +
\sqrt{p} \bigr)^2 , \qquad
\sigma_{p} = \bigl(\sqrt{n-1} + \sqrt{p} \bigr)
\biggl(\frac{1}{\sqrt{n-1}} +
\frac{1}{\sqrt{p}} \biggr)^{1/3} .
\end{equation}
Here $F_1$ denotes the orthogonal-Tracy--Widom law \cite{tw96}, the scaling
limit of the
largest eigenvalue in real Gaussian Wigner matrices. Slightly prior to
\cite{johnstone01}, as a by-product of his analysis on the random growth
model, Johansson \cite{johansson00} proved that the scaling limit for the
largest eigenvalue in the complex white Wishart matrix is the unitary
Tracy--Widom
law $F_2$. Recently, El Karoui \cite{nek-inf} extended the asymptotic regime
\eqref{eq:asymptotics} to include the cases where $n/p\to0$ or
$\infty$. For
the smallest eigenvalue, when $\gamma> 1$, Baker \textit{et al.} \cite
{bfp98} showed
that the reflection of $F_2$ about the origin is the
scaling limit for complex Wishart matrices, and Paul \cite{paul06}
gave the
Tracy--Widom limits in the case where $n/p\to\infty$ for both complex
and real
Wishart matrices.

Although this type of asymptotic result has emerged only recently in the
statistics literature, it has already found its relevance to
applications with modern data. For instance, based on
\eqref{eq:loe-lim}, Patterson \textit{et al.} \cite{patterson} developed a formal
procedure for testing the presence of population heterogeneity with SNP
(single nucleotide polymorphism) data.

From a statistical point of view, to inform the use of any asymptotic
result in practice, we need to understand how closely the asymptotic
limit approximates the finite sample distributions. In the motivating
example, this dictates the accuracy of the nominal $p$-value.

In this paper, we first establish a rate of convergence result for
the Tracy--Widom approximation to the distribution of the rescaled
largest eigenvalue, but with more carefully chosen constants than
\eqref{eq:orig-center-scale}. Set $a\wedge b = \min(a,b)$ and $m_\pm=
m\pm\frac{1}{2}$. We show that modifying the centering and scaling
constants to
%
%
\begin{equation}
\label{eq:new-center-scale}
\mu_{n,p} = \bigl(\sqrt{n_-} + \sqrt{p_-} \bigr)^2,
\qquad
\sigma_{n,p} = \bigl(\sqrt{n_-} +
\sqrt{p_-} \bigr)\biggl (\frac{1}{\sqrt{n_-}} +
\frac{1}{\sqrt{p_-}} \biggr)^{1/3}
\end{equation}
results in better approximation. The difference between the
distribution of\vspace*{2pt}
$(\lambda_1 - \mu_{n,p})/\sigma_{n,p}$ and $F_1$ reduces to the
`second order',
being $\Oh{(n\wedge p)^{-2/3}}$ rather than $\Oh{(n\wedge p)^{-1/3}}$,
that would
apply by using \eqref{eq:orig-center-scale}. See Theorem \ref{thm:main}.
Numerical work in Section~\ref{sec:np-large} suggests that the
improvement is
substantial.

Further assuming $\gamma> 1$ in \eqref{eq:asymptotics}, we find that,
with a log transform, the scaling limit of $\log\lambda_p$
is the reflected Tracy--Widom law $G_1$ (defined by $G_1(s) = 1 -
F_1(-s)$) \cite{paul06}. Moreover, with appropriate
rescaling constants, the accuracy of the limit also reaches the second
order: $\Oh{p^{-2/3}}$. See Theorem \ref{thm:small} and Section
\ref{sec:np-small}.

In the literature, El Karoui \cite{nek06} established a parallel result
for Johansson's theorem for the largest eigenvalue on the complex domain
and Choup \cite{choup06} studied the same problem via an Edgeworth
expansion approach. Recently, Johnstone \cite{johnstone07} obtained both
scaling limit and convergence rates for the extreme eigenvalues of an
$F$-matrix, on both complex and real domains. As is usual in the Random
Matrix Theory literature, results on the real domain are founded in part
on those for complex data but require significant additional constructs
and arguments. This is explained for our setting in Sections
\ref{sec:proof-thm1} and \ref{sec:proof-thm2}.

The rest of the paper is organized as follows. In Section \ref{sec:main-res},
we present theorems for both the largest and the smallest eigenvalues,
together with supporting numerical results, related statistical
settings, a
real data example and a brief discussion. Section \ref{sec:proof-thm1} proves
the theorem on the largest eigenvalue and Section \ref{sec:proof-thm2}
sketches the proof of the one on the smallest eigenvalue. Finally, Section~\ref{sec:lg} establishes necessary Laguerre polynomial asymptotics,
which is
first used without proof in Section \ref{sec:proof-thm1}. Technical details
are collected in the \hyperref[appm]{Appendix}.

\section{Main results and their applications}
\label{sec:main-res}

In this section, we first state two main theorems of this paper, which are
concerned with the convergence rates of the largest and the smallest
eigenvalues in finite Wishart matrices to their Tracy--Widom limits. The
theorems are then complemented and further justified by a series of numerical
experiments, in which the Tracy--Widom approximation is reasonably good even
when $n$ and/or $p$ are as small as $2$. After that, we review several related
statistical settings and consider a real data example. Finally, we end the
section with a brief discussion.

\subsection{Main theorems}
\label{sec:main-thm}

We begin with the largest eigenvalue, for which we have the following
rate of convergence result.
\begin{theorem}
\label{thm:main}
Let $A\sim W_p(I,n)$ with $n\neq p$ and $\lambda_1$ its largest
eigenvalue. Define~$(\mu_{n,p}, \sigma_{n,p})$ as in
\eqref{eq:new-center-scale}. Under condition\vadjust{\goodbreak} \eqref{eq:asymptotics}, for
any given $s_0$, there exists an integer~$N_0(s_0,\gamma)$, such that
when $n\wedge p\geq N_0(s_0,\gamma)$ and is even, for all $s\geq s_0$,
\[
|P\{\lambda_1\leq\mu_{n,p} + \sigma_{n,p}s\} - F_1(s) |
\leq
C(s_0)(n\wedge p)^{-2/3}\exp(-s/2) ,
\]
where $C(\cdot)$ is continuous and non-increasing.
\end{theorem}

We also obtain an analogous result for the smallest eigenvalue. Refine
condition \eqref{eq:asymptotics} to
%
%
\begin{equation}
\label{eq:asymp-1}
p \to\infty, \qquad  p + 1\leq n = n(p) \to\infty \quad \mbox{and} \quad  n/p\to
\gamma\in(1,\infty).
\end{equation}
%
Define $\mu_{n,p}^- = (\sqrt{n_-} - \sqrt{p_-} )^2$,
$\sigma_{n,p}^- = (\sqrt{n_-} -
\sqrt{p_-} ) ({1}/{\sqrt{p_-}} -
{1}/{\sqrt{n_-}} )^{1/3}$, and let
%
%
\begin{equation}
\label{eq:small-center-scale}
\tau^-_{n,p} = \frac{\sigma^-_{n,p}}{\mu^-_{n,p}}, \qquad\nu
^-_{n,p} =
\log(\mu^-_{n,p}) + \frac{1}{8} (\tau^-_{n,p} )^2.
\end{equation}
Then we have the following theorem.
\begin{theorem}
\label{thm:small}
Let $A\sim W_p(I,n)$ with $n-1\geq p$ and $\lambda_p$ as its smallest
eigenvalue. Define~
$(\nu^-_{n,p}, \tau^-_{n,p})$ as in \eqref{eq:small-center-scale}.
Under condition \eqref{eq:asymp-1}, we have
\[
\frac{\log\lambda_p - \nu^-_{n,p}}{\tau^-_{n,p}} \convd G_1
\]
with $G_1(s) = 1 - F_1(-s),$ the reflected Tracy--Widom law.

In addition, for any given $s_0$, there exists an integer
$N_0(s_0,\gamma)$, such that when $p\geq N_0(s_0,\gamma)$ and
is even, for all $s\geq s_0$,
\[
|P\{\log\lambda_p\leq\nu^-_{n,p} - \tau^-_{n,p}s\} -
G_1(-s) | \leq C(s_0) p^{-2/3}\exp(-s/2) ,
\]
where $C(\cdot)$ is continuous and non-increasing.
\end{theorem}

%
\begin{table*}
\tabcolsep=0pt
\caption{Simulations for finite $n\times p$ vs. Tracy--Widom limit: the
largest eigenvalue.
For each $(n,p)$ combination,
we show in the first row empirical cumulative probabilities for
$\lambda_1$, rescaled by \protect\eqref{eq:new-center-scale}, and the
second row, with parentheses, rescaled by
\protect\eqref{eq:orig-center-scale}, both computed from $R = 40\,000$
repeated draws from $W_p(n, I)$ using the method in \protect\cite{de02}. Conventional
significance levels are highlighted in bold font and
the last row gives approximate standard errors based on binomial
sampling. $F_1$ was computed by the method in \protect\cite{persson05} with
percentiles obtained via inverse interpolation}\label{table:prob-tb}
\begin{tabular*}{\textwidth}{@{\extracolsep{\fill}}ld{2.4}d{2.4}d{2.4}d{2.4}d{2.4}d{2.4}d{2.4}d{2.4}l@{}}
\hline
Percentiles & \multicolumn{1}{l}{$-3.8954$} & \multicolumn{1}{l}{$-3.1804$} & \multicolumn{1}{l}{$-2.7824$} & \multicolumn{1}{l}{$-1.9104$}
& \multicolumn{1}{l}{$-1.2686$}
& \multicolumn{1}{l}{$-0.5923$} & \multicolumn{1}{l}{$0.4501$} & \multicolumn{1}{l}{$0.9793$} &
\multicolumn{1}{l@{}}{$2.0234$}\\
\hline
TW & 0.01 & 0.05 & 0.10 & 0.30 & 0.50 & 0.70 & \multicolumn{1}{l}{\hphantom{0}\textbf{0.90}} & \multicolumn{1}{l}{\hphantom{0}\textbf{0.95}}
& \multicolumn{1}{l@{}}{\hphantom{(}\textbf{0.99}} \\
[5pt]
$2\times2$ & 0.000 & 0.000 & 0.000 & 0.034 & 0.379 & 0.690 & \multicolumn{1}{l}{\hphantom{0}\textbf{0.908}}
& \multicolumn{1}{l}{\hphantom{0}\textbf{0.953}} & \multicolumn{1}{l@{}}{\hphantom{(}\textbf{0.988}} \\
& (0.000) & (0.000) & (0.000) & (0.015) & (0.345) & (0.669) & (0.902) & (0.950)
& (0.987) \\
$5\times5$ & 0.000 & 0.002 & 0.021 & 0.218 & 0.465 & 0.702 & \multicolumn{1}{l}{\hphantom{0}\textbf{0.908}}
& \multicolumn{1}{l}{\hphantom{0}\textbf{0.954}} & \multicolumn{1}{l@{}}{\hphantom{(}\textbf{0.989}} \\
& (0.000) & (0.002) & (0.020) & (0.213) & (0.460) & (0.698) & (0.907) & (0.953)
& (0.989) \\
$25\times25$ & 0.003 & 0.031 & 0.075 & 0.280 & 0.492 & 0.700 & \multicolumn{1}{l}{\hphantom{0}\textbf{0.902}} & \multicolumn{1}{l}{\hphantom{0}\textbf{0.951}} & \multicolumn{1}{l@{}}{\hphantom{(}\textbf{0.990}} \\
& (0.003) & (0.030) & (0.075) & (0.280) & (0.491) & (0.699) & (0.902) & (0.951)
& (0.990) \\
$100\times100$ & 0.007 & 0.041 & 0.091 & 0.294 & 0.501 & 0.704 & \multicolumn{1}{l}{\hphantom{0}\textbf{0.902}} & \multicolumn{1}{l}{\hphantom{0}\textbf{0.951}} & \multicolumn{1}{l@{}}{\hphantom{(}\textbf{0.990}} \\
& (0.007) & (0.041) & (0.091) & (0.294) & (0.501) & (0.704) & (0.902) & (0.951)
& (0.990) \\
[5pt]
$8\times2$ & 0.000 & 0.001 & 0.012 & 0.196 & 0.456 & 0.702 & \multicolumn{1}{l}{\hphantom{0}\textbf{0.909}}
& \multicolumn{1}{l}{\hphantom{0}\textbf{0.955}} & \multicolumn{1}{l@{}}{\hphantom{(}\textbf{0.990}} \\
& (0.000) & (0.004) & (0.031) & (0.270) & (0.532) & (0.754) & (0.928) & (0.964)
& (0.992) \\
$20\times5$ & 0.001 & 0.018 & 0.054 & 0.259 & 0.483 & 0.704 & \multicolumn{1}{l}{\hphantom{0}\textbf{0.906}}
& \multicolumn{1}{l}{\hphantom{0}\textbf{0.954}} & \multicolumn{1}{l@{}}{\hphantom{(}\textbf{0.990}} \\
& (0.002) & (0.028) & (0.073) & (0.303) & (0.531) & (0.737) & (0.921) & (0.962)
& (0.992) \\
$100\times25$ & 0.006 & 0.040 & 0.088 & 0.292 & 0.498 & 0.701 & \multicolumn{1}{l}{\hphantom{0}\textbf{0.901}} & \multicolumn{1}{l}{\hphantom{0}\textbf{0.950}} & \multicolumn{1}{l@{}}{\hphantom{(}\textbf{0.989}} \\
& (0.008) & (0.047) & (0.100) & (0.314) & (0.523) & (0.721) & (0.910) & (0.955)
& (0.991) \\
$400\times100$ & 0.009 & 0.048 & 0.096 & 0.299 & 0.502 & 0.702 & \multicolumn{1}{l}{\hphantom{0}\textbf{0.902}} & \multicolumn{1}{l}{\hphantom{0}\textbf{0.951}} & \multicolumn{1}{l@{}}{\hphantom{(}\textbf{0.990}} \\
& (0.010) & (0.053) & (0.104) & (0.312) & (0.516) & (0.714) & (0.908) & (0.954)
& (0.991) \\
[5pt]
$500\times5$ & 0.010 & 0.049 & 0.098 & 0.296 & 0.502 & 0.705 & \multicolumn{1}{l}{\hphantom{0}\textbf{0.906}} & \multicolumn{1}{l}{\hphantom{0}\textbf{0.955}} & \multicolumn{1}{l@{}}{\hphantom{(}\textbf{0.990}} \\
& (0.020) & (0.083) & (0.150) & (0.385) & (0.589) & (0.772) & (0.933) & (0.969)
& (0.994) \\
$1000\times10$ & 0.010 & 0.051 & 0.101 & 0.300 & 0.504 & 0.707 & \multicolumn{1}{l}{\hphantom{0}\textbf{0.902}} & \multicolumn{1}{l}{\hphantom{0}\textbf{0.952}} & \multicolumn{1}{l@{}}{\hphantom{(}\textbf{0.991}} \\
& (0.017) & (0.077) & (0.138) & (0.366) & (0.571) &
(0.757) & (0.923) & (0.963) & (0.994) \\
$5000\times5$ & 0.012 & 0.056 & 0.107 & 0.307 & 0.509 & 0.707 & \multicolumn{1}{l}{\hphantom{0}\textbf{0.905}} & \multicolumn{1}{l}{\hphantom{0}\textbf{0.953}} & \multicolumn{1}{l@{}}{\hphantom{(}\textbf{0.992}} \\
& (0.027) & (0.097) & (0.169) & (0.402) & (0.602) & (0.779) & (0.933) & (0.969)
& (0.994) \\
$10\,000\times10$ & 0.012 & 0.055 & 0.108 & 0.308 & 0.504 & 0.706 & \multicolumn{1}{l}{\hphantom{0}\textbf{0.905}} & \multicolumn{1}{l}{\hphantom{0}\textbf{0.954}} & \multicolumn{1}{l@{}}{\hphantom{(}\textbf{0.991}} \\
& (0.021) & (0.084) & (0.150) & (0.378) & (0.580) & (0.763) & (0.929) & (0.967)
& (0.994) \\
[5pt]
$2\times\mathrm{SE}$ & 0.001 & 0.002 & 0.003 & 0.005 & 0.005 & 0.005 &0.003 & 0.002 & \hphantom{(}0.001\\
\hline
\end{tabular*}
\end{table*}

While we only prove rigorous bounds for even $p$, numerical experiments
show that the approximation works just as well in the odd case, and for
the largest eigenvalue, also in the square case. See Tables
\ref{table:prob-tb} and \ref{table:prob-tb-small}.


\subsection{Numerical performance}
\label{sec:num-perf}

An important motivation for the current study is to promote practical
use of the Tracy--Widom approximation. To this end, we conduct here a set
of experiments to investigate its numerical quality.

\subsubsection{The largest eigenvalue}
\label{sec:np-large}

\paragraph*{Distributional approximation}
We first computed the empirical cumulative probabilities of $\lambda
_1$ (after
rescaling), at a collection of $F_1$ percentiles, using $R = 40\,000$
replications. This\vadjust{\goodbreak} is done for three different categories of ($n,p$)
combinations: (1) the square case, where $n = p = 2$, $5$, $25$ and
$100$; (2)
the rectangular case, where $p = 2$, $5$, $25$ and $100$ and $n/p$ is fixed
at 4:1; (3) the `thin' case, where $p = 5$ and $10$ but $n/p$ could be
as high
as 100:1 and 1000:1{}. In some sense, this category could also be
thought of as in the situation where $n/p\to\infty$ as discussed in
\cite{nek-inf}. For comparison purpose, we rescaled $\lambda_1$
using both
the new constants \eqref{eq:new-center-scale} and the old ones
\eqref{eq:orig-center-scale}. The results are summarized in Table
\ref{table:prob-tb}.

\begin{table*}
\tabcolsep=0pt
\caption{Simulations for finite
$n\times p$ vs. Tracy--Widom
limit: the smallest eigenvalue. For each $(n,p)$ combination,
empirical cumulative probabilities are computed for
$(\log\lambda_p - \nu^-_{n,p})/\tau^-_{n,p}$ using $R = 40\,000$ draws
from $W_p(I, n)$. Methods for sampling, computing $F_1$ and
obtaining percentiles are the same as in Table
\protect\ref{table:prob-tb}. Conventional significance levels are
highlighted in bold font and the last line gives approximate
standard errors based on binomial sampling}\label{table:prob-tb-small}
\begin{tabular*}{\textwidth}{@{\extracolsep{\fill}}ld{1.3}d{1.3}d{1.3}d{1.3}d{1.3}d{1.3}lll@{}}
\hline
Percentiles & \multicolumn{1}{l}{3.8954} & \multicolumn{1}{l}{3.1804} & \multicolumn{1}{l}{2.7824} & \multicolumn{1}{l}{1.9104}
& \multicolumn{1}{l}{1.2686} & \multicolumn{1}{l}{0.5923} &
\multicolumn{1}{l}{$-$0.4501} & \multicolumn{1}{l}{$-$0.9793} &
\multicolumn{1}{l@{}}{$-$2.0234}\\
\hline
RTW & 0.99 & 0.95 & 0.90 & 0.70 & 0.50 & 0.30 & \textbf{0.10} & \textbf{0.05}
& \textbf{0.01} \\
[5pt]
$4\times2$ & 1.000 & 1.000 & 0.998 & 0.893 & 0.625 & 0.326 & \textbf
{0.087} & \textbf{0.041} & \textbf{0.009} \\
$10\times5$ &0.999 & 0.995 & 0.976 & 0.798 & 0.555 & 0.310 & \textbf{0.095}
& \textbf{0.047} & \textbf{0.011} \\
$50\times25$ & 0.997 & 0.973 & 0.931 & 0.728 & 0.515 & 0.302 & \textbf
{0.097} & \textbf{0.048} & \textbf{0.010} \\
$200\times100$ & 0.993 & 0.960 & 0.913 & 0.713 & 0.509 & 0.306 & \textbf
{0.103} & \textbf{0.050} & \textbf{0.010} \\
[5pt]
$8\times2 $ & 1.000 & 0.992 & 0.969 & 0.792 & 0.554 & 0.314 & \textbf
{0.095} & \textbf{0.046} & \textbf{0.010}\\
$20\times5$ & 0.999 & 0.977 & 0.939 & 0.740 & 0.522 & 0.301 & \textbf{0.096}
& \textbf{0.047} & \textbf{0.009} \\
$100\times25$ & 0.993 & 0.960 & 0.915 & 0.713 & 0.505 & 0.298 & \textbf
{0.098} & \textbf{0.048} & \textbf{0.009} \\
$400\times100$ & 0.992 & 0.954 & 0.904 & 0.701 & 0.500 & 0.298 & \textbf
{0.100} & \textbf{0.049} & \textbf{0.010} \\
[3pt]
$2\times\mathrm{SE}$ & 0.001 & 0.002 & 0.003 & 0.005 & 0.005 & 0.005 &0.003 & 0.002 & 0.001\\
\hline
\end{tabular*}
\end{table*}

Numerical accuracy with the new constants could be viewed from two
aspects. First, for the conventional significance levels of $10\%$,
$5\%$ and $1\%$ that correspond to right tails of the distributions,
the approximation looks good even when $p$ is as small as $2$! In
addition, it improves as $p$ becomes larger and starts to match the
finite distributions almost exactly when $p$ is no greater than $25$. See
the last three columns of Table \ref{table:prob-tb}. Second, when $p$ is
large, for instance, in the $100\times100$ and $400\times100$ cases,
$F_1$ provides reasonable approximation over the whole range
of interest.

As regards the comparison between different rescaling constants, neither
choice seems superior to the other in the square cases (see the first
block of Table \ref{table:prob-tb}). However, when the ratio $n/p$ is
changed to 4:1 or higher (see the second and the third blocks), the
improvement by using new constants \eqref{eq:new-center-scale} is
self-evident.

As a remark, better performance on right tails and improvement by using
the new constants, as reflected in this simulation study, agree well
with the mathematical statement in Theorem \ref{thm:main}.

\paragraph*{Approximate percentiles}
We can also use $F_1$ to calculate approximate percentiles for the
finite distributions, whose accuracy can be measured by the relative
error \mbox{$r_{\aa} = \theta^{TW}_{\aa}/\theta_{\aa}-1$}. Here, $\theta
_{\aa}$
is the exact $100\aa$th percentile of the rescaled largest eigenvalue
in the finite $n\times p$ model and $\theta^{TW}_{\aa}$ is its
counterpart from $F_1$.

%
\begin{figure}
\centering
\begin{tabular}{@{}cc@{}}

\includegraphics{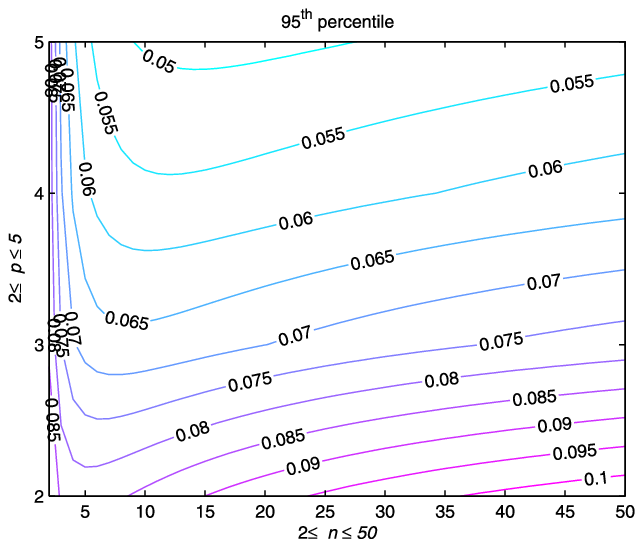}
&{\includegraphics{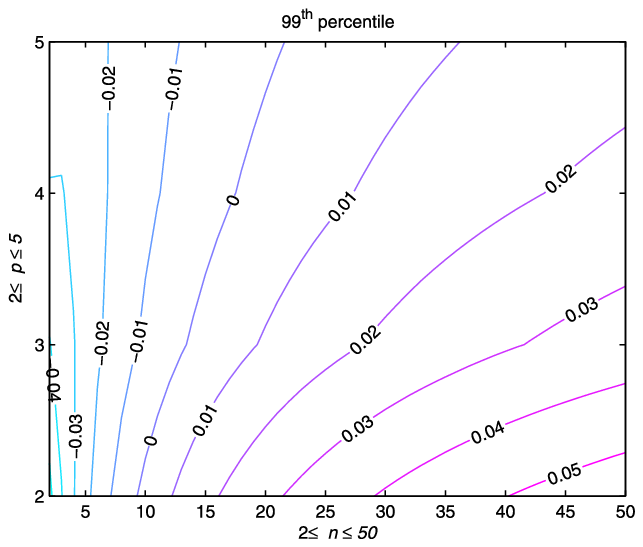}}\\
(a)&(b)
\end{tabular}
\vspace*{-3pt}
\caption{Plots of relative errors $r_\aa$ for
approximate percentiles using $F_1$: \textup{(a)} $95$th percentile; \textup{(b)}
$99$th percentile. Exact finite $n\times p$ distributions are
computed in MATLAB using Koev's implementation \protect\cite{koev} and
$F_1$ is computed using the method in \protect\cite{persson05}. The
percentiles are obtained from inverse interpolation.}\label{fig:rr}
\vspace*{-3pt}
\end{figure}

In Figure~\ref{fig:rr}, we plot $r_\aa$ for $\aa= 0.95$ and $0.99$,
with $p$
ranging from $2$ to $5$ and $n$ from~$2$ to $50$. Although $n \wedge p$
is no
greater than $5$, the approximation is reasonably satisfactory. For the
$95$th percentile, $|r_{0.95}|$ ranges from $5\%$ to $10\%$ for most cases
and slightly exceeds $10\%$ only when $p = 2$ and the $n/p$ ratio is
high. The
approximation works even better for the $99$th percentile, with $|
r_{0.99}| \leq5\%$ for most cases. Due to computational limitation~\cite{koev}, we could not obtain exact percentiles when $n$ and $p$
are large.
We expect the approximate percentiles to become more accurate as a consequence
of better distributional approximation.

\subsubsection{The smallest eigenvalue}
\label{sec:np-small}

For the smallest eigenvalue, we perform a simulation study to
investigate the
distributional approximation. We chose two $n/p$ ratios: 2:1 and 4:1, both
with $p = 2$, $5$, $25$ and~$100$. For each $(n,p)$ combination, we
used $R =
40\,000$ replications. The simulation results shown in Table
\ref{table:prob-tb-small} demonstrate similar performance as in the
case of
the largest eigenvalue and agree well with Theorem \ref{thm:small}.

\subsection{Related statistical settings}
\label{sec:stat-set}

Here, we review several settings in multivariate statistics to which our
results are applicable. Throughout the subsection, we only use the
largest eigenvalue to illustrate.

\subsubsection*{Principal component analysis}
Suppose that $X =
[{X}_1,\ldots, {X}_n]'$ is a Gaussian data matrix. Write the sample
covariance matrix ${S} = (n-1)^{-1}{X}'{H}{X},$ where ${H} =
{I}-n^{-1}\mathbf{1}\mathbf{1}'$ is the centering matrix and principal
component analysis (PCA) looks for a sequence of standardized vectors
${a}_1, \ldots, {a}_p$ in $\mathbb{R}^p$, such that ${a}_i$ successively
solves the following optimization problem:
\[
\max\{{a}'{S}{a}\dvtx {a}'{a}_j =
0, j\leq i \} ,
\]
where ${a}_0$ is the zero vector. Then, successive sample eigenvalues
$\hat\ell_1\geq\cdots\geq\hat\ell_p$ satisfy $\hat\ell_i =
{a}'_i{S}{a}_i$.

One basic question in PCA application is testing the hypothesis of isotropic
variation, that is, the population\vadjust{\goodbreak} covariance matrix $\Sigma= \tau^2
I$. For
simplicity, assume that $\tau^2 = 1$ (otherwise we divide $S$ by $\tau^2$).
Then $(n-1){S}\sim W_p({I}, n-1)$. The largest eigenvalue $\hat\ell
_1$ of
${S}$ is a natural test statistic under the union intersection
principle. Our
result applies to $(n-1)\hat\ell_1$. If $\tau^2$ is unknown, we
could estimate
it by $\trace{S}/p$. See \cite{nadler10}.

\subsubsection*{Testing that a covariance matrix equals a specified matrix}
Suppose that $X = [X_1,\ldots,X_n]'$ has as its row vectors i.i.d. samples
from the $N_p(\mu,\Sigma)$ distribution. We want to test the hypothesis
$H_0\dvt \Sigma= \Sigma_0$, where $\Sigma_0$ is a specified positive
definite matrix.


Suppose $\mathbf{\mu}$ is unknown, and let ${S} = (n-1)^{-1}{X}'{H}{X}$
be the sample covariance matrix. The union intersection test uses the
largest eigenvalue of ${\Sigma}_0^{-1}{S}$, denoted by
$\hat\ell_1({\Sigma}_0^{-1}{S})$, as the test statistic
\cite{mkb}, page~130. Observe that $\hat\ell_1({\Sigma}_0^{-1}{S}) =
\hat\ell_1({\Sigma}_0^{-1/2}{S}{\Sigma}_0^{-1/2})$. Under $H_0$,
$(n-1){\Sigma}_0^{-1/2}{S}{\Sigma}_0^{-1/2} \sim W_p({I}, n-1)$. So,
our result is available for $(n-1)\hat\ell_1({\Sigma}_0^{-1}{S})$.

\subsubsection*{Singular value decomposition} For $X$ a real $n\times p$
matrix, there exist orthogonal matrices $U(n\times n)$ and
$V(p\times p)$, such that
\[
X = UDV^T,
\]
where $D = \operatorname{diag}(d_1, \ldots, d_{n\wedge p}) \in
\mathbb{R}^{n\times p}$ and $d_1\geq\cdots\geq d_{n\wedge p}\geq
0$. This representation is called the singular value decomposition of
$X$ \cite{hj85}, Theorem 7.3.5, with $d_i$ the $i$th singular value of
$X$. Theorem \ref{thm:main} then provides an accurate distributional
approximation for $d_1^2$ when the entries of $X$ are independent
standard normal random variables.

\subsection{The score data example}
\label{sec:data}

We consider now the score data example extracted from \cite{mkb}. The
data set consists of the scores of 88 students on 5 subjects
(mechanics, vectors, algebra, analysis and statistics). Taking account
of centering, we have $n = 87$ and $p = 5$.

One might expect that there are several common factors that determine
the students' performance on the tests. Moreover, one might assume that
the joint effects of the common factors are observed in isotropic
noises, in which case the covariance structure of the scores (after
proper diagonalization) follows a spiked model $\Sigma=
\tau^2\Sigma_m$, where $\tau^2 > 0$, $\Sigma_m =
\operatorname{diag}(\ell_1,\ldots,\ell_m,1,\ldots,1)$ and $0\leq
m\leq
4$. (Note that the model $\Sigma= \tau^2\Sigma_4$ is the saturated
model and is indistinguishable from $\Sigma= \tau^2\Sigma_5$.) To
determine $m$, we are led to test a nested sequence of hypotheses $H_k\dvt
\Sigma= \tau^2\Sigma_m$ with some $m\leq k$, for $0\leq k\leq3$.

To compute the $p$-value of testing $H_k$, we could (i) estimate
$\tau^2$ by $\hat\tau^2_{p-k}$\vspace*{-2pt} as the mean of the $p - k$ smallest
sample eigenvalues; (ii) construct the test statistic\vspace*{2pt} as $T_k =
(n\hat\ell_{k+1}/\hat\sigma^2_{p-k} - \mu_{n,p-k})/\sigma_{n,p-k}$;
(iii) report $F_1(T_k)$ as the approximate conservative $p$-value. Step
(iii) is justified as follows. Let $\mathcal{L}(\lambda_j|n, p,
\Sigma)$
denote the law of the $j$th largest sample eigenvalue of a $W_{p}(n,
\Sigma)$ matrix. By the interlacing properties of the eigenvalues
\cite{hj85}, Theorem 7.3.9 (see also \cite{johnstone01}, Proposition
1.2), $\mathcal{L}(\lambda_1|n,p - m,I_{p-m})$ could be
used to compute the conservative $p$-value for the null distribution
$\mathcal{L}(\lambda_{k+1}|n,p,\Sigma_m)$ for all $k\geq m$, which is
further approximated by $F_1$. We summarize the values of $T_k$ and the
corresponding $p$-values in Table \ref{table:score}.

%
\begin{table*}
\tabcolsep=0pt
\tablewidth=308pt
\caption{The test statistics $T_k$ and the corresponding $p$-values
$F_1(T_k)$ calculated using new centering and scaling constants
(\protect\ref{eq:new-center-scale}) and old constants
(\protect\ref{eq:orig-center-scale}) for the score data}\label{table:score}
\begin{tabular*}{308pt}{@{\extracolsep{\fill}}lllll@{}}
\hline
& $H_0$ & $H_1$ & $H_2$ & $H_3$ \\
\hline
$T_k$ (new) & 14.5934 & 4.3162 & 0.4535 & 1.4949 \\
$p$-value (new) & $<10^{-6}$ & $1.1\times10^{-4}$ & 0.0996 & 0.0235\\
 [5pt]
$T_k$ (old) & 14.4740 & 4.1155 & 0.1803 & 1.1897 \\
$p$-value (old) & $<10^{-6}$ & $1.7\times10^{-4}$ & 0.1376 & 0.0371\\
\hline
\end{tabular*}
\end{table*}

From Table \ref{table:score}, we could see a noticeable difference
between the values of $T_k$ and the corresponding $p$-values by using
different rescaling constants. The $p$-values obtained from the new
constants are typically smaller than those from the old
constants. Noting that the $p$-values are already conservative, the new
constants \eqref{eq:new-center-scale} prevent further unnecessary
conservativeness that would otherwise be caused by the old constants in
this example.

\subsection{Discussion}
\label{sec:discuss}

We discuss below two issues related to our results.

\subsubsection*{Log transform} One notable difference between
Theorems \ref{thm:main} and \ref{thm:small} is the logarithmic
transformation of the smallest eigenvalue before scaling.

Indeed, for the largest eigenvalue, a similar $\Oh{N^{-2/3}}$
convergence rate
can be obtained for the distribution of $(\log\lambda_1 -
\nu_{n,p})/\tau_{n,p}$, with $\nu_{n,p} = \log(\mu_{n,p})$ and
$\tau_{n,p} =
\sigma_{n,p}/\mu_{n,p}$. However, when $n$ or $p$ is small, its numerical
results are not as good as those obtained from direct scaling. In comparison,
for the smallest eigenvalue, the transform yields substantial numerical
improvement. Therefore, we recommend the log transform for the smallest
eigenvalue.

As no theoretical analysis justifying the choice of the transform is currently
available, we attempt some heuristics in the following. First, observe that
sample covariance matrices are positive semidefinite. So, for $\lambda_p$,
the hard lower bound at $0$ truncates the left tail of its density
function on
any linear scale, and hence obstructs the asymptotic approximation by $G_1$
that is supported on the whole real line. However, by a map $x\mapsto
\log{x}$, we map the support to the whole real line and avoid the `hard
edge' effect. The largest eigenvalue does not necessarily benefit from this
transform, for it is on the `soft edge', that is, the right edge of the
covariance matrix spectrum, which does not have a deterministic upper
bound. Such heuristics are supported by related studies on Gaussian Wigner
matrices \cite{jm09} and $F$-matrices \cite{johnstone07}.

\subsubsection*{Software} There have been works on the numerical
evaluation of the
Tracy--Widom distributions~\cite{bejan,persson05,bornemann11} and
the exact
finite $n\times p$ distributions of the extreme eigenvalues \cite
{koev-web,koev}.
In addition, the author and colleagues have developed an \texttt{R}
package \texttt{RMTstat} \cite{roylab} that is intended to provide an
interface for using the Tracy--Widom approximation in multivariate statistical
analysis.

\vspace*{-2pt}\section{The largest eigenvalue}\vspace*{-2pt}
\label{sec:proof-thm1}

This section is devoted to the proof of Theorem \ref{thm:main}. We use
the operator norm convergence framework developed in \cite{tw05}, for
the joint eigenvalue distribution of white Wishart matrices is
essentially the same as the Laguerre orthogonal ensemble in random
matrix theory (RMT).

In the proof, we first give the determinantal representations for the
finite and limiting distribution functions and work out explicit
formulas for related kernels, in which Widom's formula
\eqref{eq:widom-S1} plays the central role. Then, a Lipschitz-type
inequality shows that the difference in determinants is bounded by the
difference in kernels. The representation of the finite sample kernel
involves weighted generalized Laguerre polynomials, while that of the
limiting kernel uses Airy function. A decomposition of the
kernel difference then enables us to transfer bounds on the convergence of
Laguerre polynomials to Airy function to bounds on the kernel
difference and eventually to bounds on the difference of the
probabilities.


\subsection{Determinantal laws}
\label{sec:det-law}

Following RMT notational convention, we replace the dimension parameter
$p$ of a white Wishart matrix $A$ by $N$, and use $x_i$ instead of
$\lambda_i$ to denote its eigenvalues. Henceforth, we assume that $N$ is
even, $n=n(N) \geq N+1$ and $n/N\to\gamma\in[1,\infty)$ as $N\to
\infty$. The cases $\gamma\in(0,1]$ are easily obtained by
interchanging $n$ and $N$.

In the RMT literature, for an integer $N \geq2$ and any $\aa> -1$, the
Laguerre orthogonal ensemble with parameters $N$ and $\aa$, denoted by
$\operatorname{LOE}(N,\aa)$, refer to joint eigenvalue density
%
%
\begin{equation}
\label{eq:loe-pdf}
\tilde{p}_N(x_1,\ldots, x_N) =
\frac{1}{{Z}_{N,\aa}}\prod_{1\leq j<k\leq N}(x_j - x_k)
\prod_{j=1}^N x_j^{\aa}\mathrm{e}^{-x_j/2},
\end{equation}
where $x_1\geq\cdots\geq x_N\geq0$. If further $\aa$ is a
non-negative integer, \eqref{eq:loe-pdf} matches the density function
of ordered eigenvalues $x_1\geq\cdots\geq x_N \geq0$ from a white
Wishart matrix $A\sim W_N(I,n)$, with
%
%
\begin{equation}
\label{eq:ch-var}
\aa= n - N - 1.
\end{equation}
Henceforth, we identify the LOE($N,\aa$) model with eigenvalues of
$A\sim
W_N(I,n)$ by~\eqref{eq:ch-var}. Thinking of $\aa$ and $n$ as functions
of $N$, in what follows we sometimes drop explicit dependence of
certain quantities on them.\vadjust{\goodbreak}

For LOE($N,\aa$), \cite{tw98}, Section
9, features the following determinantal formula
%
%
\begin{equation}
\label{eq:F-N-1}
\tilde{F}_{N,1}(x') = P\{x_1\leq x'\} = \sqrt{\det(I-K_N\chi)}.
\end{equation}
Here $\chi= \id_{x>x'}$ and $K_N$ is an operator with $2\times2$
matrix kernel
%
%
\begin{equation}
\label{eq:K-N-kernel}
K_{N}(x,y) = (L S_{N,1})(x,y) + K^\varepsilon(x,y),
\end{equation}
where
\[
L =
\pmatrix{ I & -\partial_2 \cr \varepsilon_1 & T
}
,\qquad K^\varepsilon=
\pmatrix{ 0 & 0 \cr -\varepsilon(x-y) & 0
}
.
\]
In $L$, $\partial_2$ is the differential operator with respect to the
second argument, $\varepsilon_1$ is the convolution operator acting on the
first argument with the kernel $\varepsilon(x-y)=\frac
{1}{2}\operatorname{sgn}(x-y)$
and $TK(x,y) = K(y,x)$ for any kernel $K$.

To give an explicit formula for $S_{N,1}$, introduce the generalized
Laguerre polynomials $\{L_k^\aa\}_{k=0}^\infty$ (\cite{szego},
Chapter~V), which are orthogonal on $[0,\infty)$ with weight function
$x^\aa
\mathrm{e}^{-x}$. The normalized and weighted versions of them become
%
%
\begin{equation}
\label{eq:phi-k}
\phi_k(x;\aa) =
h_k^{-1/2}x^{\aa/2}\mathrm{e}^{-x/2}L_k^{\aa}(x), \qquad k = 0,\ldots,
\end{equation}
with $h_k = \int_0^\infty L_k^\aa(x)^2x^\aa \mathrm{e}^{-x}\,\mathrm{d}x =
(k+\aa)!/k!$. Widom \cite{widom99} derived a formula for $S_{N,1}$,
which can be rewritten in a form more convenient to us
\cite{adler00}, equation (4.3), as
%
%
\begin{eqnarray}
\label{eq:widom-S1}
S_{N,1}(x,y) &=& S_{N,2}(x,y) +   \frac{N!}{4\Gamma(N+{\aa})}
x^{{\aa}/2} \mathrm{e}^{-x/2} \biggl[\frac{\mathrm{d}}{\mathrm{d}x} L_N^{{\aa}}(x)
\biggr]\nonumber
\\[-8pt]
\\[-8pt]
&&\hphantom{S_{N,2}(x,y) +}
{}\times\int_0^\infty
\operatorname{sgn}(y-z) z^{{\aa}/2-1} \mathrm{e}^{-z/2} [L_{N}^{{\aa}}(z) -
L_{N-1}^{{\aa}}(z)]\,\mathrm{d}z,
\nonumber
\end{eqnarray}
where $S_{N,2}$ is the unitary correlation kernel
\[
S_{N,2}(x,y) = \sum_{k = 0}^{N-1}\phi_k(x;\aa)\phi_k(y;\aa).
\]

Let $a_N = \sqrt{N(N+\aa)}$, and define as in \cite{nek06}, Section 2,
functions
%
%
\begin{eqnarray}
\label{eq:phi-psi}
\phi(x;\aa) & =& (-1)^N\sqrt{\frac{a_N}{2}}\phi_N(x;\aa-
1)x^{-1/2}\id_{x \geq0},\nonumber
\\[-8pt]
\\[-8pt]
\psi(x;\aa) & =&
(-1)^{N-1}\sqrt{\frac{a_N}{2}}\phi_{N-1}(x;\aa+ 1)x^{-1/2}\id_{x
\geq0}.
\nonumber
\end{eqnarray}
Write $a\diamond b$ for the operator with kernel $(a\diamond b)(x,y) =
\int_0^\infty
a(x+z)b(y+z)\,\mathrm{d}z$. 
Then $S_{N,2}$ has the integral representation \cite{johnstone01,nek06}
%
%
\begin{equation}
\label{eq:S-N-2-int}
S_{N,2}(x,y) = \int_0^\infty\phi(x+z)\psi(y+z) +
\psi(x+z)\phi(y+z)\,\mathrm{d}z = (\phi\diamond\psi+ \psi\diamond\phi)(x,y).
\end{equation}
By \cite{szego}, equations (5.1.13) and (5.1.14), the second term on
the right-hand side of
\eqref{eq:widom-S1} equals\vspace*{-2pt}
\[
-\frac{N!}{4\Gamma(N+\tilde{\aa})}x^{\aa/2}\mathrm{e}^{-x/2}L_{N-1}^{\aa+1}(x)
\int_0^\infty\operatorname{sgn}(y-z)z^{\aa/2 - 1}\mathrm{e}^{-z/2}L_N^{\aa-1}(z)\,\mathrm{d}z
= \psi(x)(\varepsilon\phi)(y).\vspace*{-2pt} 
\]
Hence, we obtain\vspace*{-2pt}
%
%
\begin{equation}
\label{eq:central}
S_{N,1}(x,y) = S_{N,2}(x,y) + \psi(x)(\varepsilon\phi)(y)\vspace*{-2pt}
\end{equation}
with $S_{N,2}(x,y)$ given in \eqref{eq:S-N-2-int}. Together with
\eqref{eq:F-N-1} and \eqref{eq:K-N-kernel}, this gives the determinantal
representation of the finite sample distribution on the original scale.

The Tracy--Widom limit has a corresponding determinantal representation
\cite{tw05}\vspace*{-2pt}
%
%
\begin{equation}
\label{eq:F-1-2}
F_1(s') = \sqrt{\det(I-K_{\mathrm{GOE}}f)},\vspace*{-2pt}
\end{equation}
where $f = \id_{s > s'}$ and the operator $K_{\mathrm{GOE}}$ has the matrix kernel\vspace*{-2pt}
\[
K_{\mathrm{GOE}}(s,t) =
\pmatrix{
S(s,t) & SD(s,t)\vspace*{-2pt}\cr
IS(s,t) & S(t,s)
}
+ K^\varepsilon(s,t).\vspace*{-2pt}
\]
Introduce the right tail integration operator $\tee$ as in
\cite{johnstone07}, where $(\tee g)(s) =\int_s^\infty g(u)\,\mathrm{d}u$ and for
kernel $K(s,t)$, $(\tee_1 K)(s,t) = \int_s^\infty K(u,t)\,\mathrm{d}u$. Also
write $a\otimes b$ for the rank one operator with kernel $(a\otimes
b)(s,t) = a(s)b(t)$. Then the entries of $K_{\mathrm{GOE}}$ are\vspace*{-2pt}
%
%
\begin{eqnarray}
\label{eq:k-goe-entries}
S(s,t) & =&\bigl (S_A - \tfrac{1}{2}\Ai\otimes\tee\Ai\bigr)(s,t) + \tfrac
{1}{2}\Ai(s),\nonumber\\[-1pt]
SD(s,t) & =& -\partial_2 \bigl(S_A(s,t) - \tfrac{1}{2}\Ai\otimes\tee\Ai
\bigr)(s,t),\\[-1pt]
IS(s,t) & =& -\tee_1\bigl(S_A -\tfrac{1}{2}\Ai\otimes\tee\Ai\bigr)(s,t) -
\tfrac{1}{2}(\tee\Ai)(s) + \tfrac{1}{2}(\tee\Ai)(t).\nonumber\vspace*{-2pt}
\end{eqnarray}
Here $S_A(s,t) = (\Ai\diamond\Ai)(s,t)$ is the
Airy kernel, and $\Ai(\cdot)$ is the Airy function (\cite{olver74},
page~53,
equation (8.01)).

Let $G = \frac{1}{\sqrt{2}}\Ai$, and define matrix operators\vspace*{-2pt}
\[
\tilde{L} =
\pmatrix{
I & -\partial_2 \vspace*{-1pt}\cr
-\tee_1 & T
}
,\qquad
L_1 =
\pmatrix{
I & 0 \vspace*{-1pt}\cr
-\tee_1 & 0
}
, \qquad
L_2 =
\pmatrix{
0 & 0 \vspace*{-1pt}\cr
\tee_2 & I
}
.\vspace*{-2pt}
\]
We can write $K_{\mathrm{GOE}}$ in a compact form as\vspace*{-2pt}
%
%
\begin{equation}
\label{eq:K-GOE}
K_{\mathrm{GOE}} = \tilde{L}(S_A - G \otimes\tee G) +
L_1\bigl(G\otimes\tfrac{1}{\sqrt{2}}\bigr) + L_2\bigl(\tfrac{1}{\sqrt{2}}\otimes
G\bigr) + K^\varepsilon.\vspace*{-3pt}
\end{equation}

\subsection{Rescaling the finite sample kernel}\vspace*{-3pt}
\label{sec:rescale}

Under the current RMT notation, the rescaling constants
\eqref{eq:new-center-scale} are translated to\vspace*{-2pt}
%
%
\begin{equation}
\label{eq:tmu-tsigma}
\mu_{n,N} = \bigl(\sqrt{n_-} + \sqrt{N_-} \bigr)^2,\qquad
\sigma_{n,N} = \bigl(\sqrt{n_-} + \sqrt{N_-} \bigr)
\biggl(\frac{1}{\sqrt{n_-}} + \frac{1}{\sqrt{N_-}} \biggr)^{1/3}.\vspace*{-2pt}
\end{equation}
Introduce the linear transformation $\tau(s) = \tauexpr{s}$ and let
$F_{N,1}(\cdot) = \tilde{F}_{N,1}(\tau(\cdot))$ be the distribution
function of $\tau^{-1}(x_1)$, that is, the largest eigenvalue of
$A\sim
W_N(I,n)$, rescaled by \eqref{eq:tmu-tsigma}.\vadjust{\goodbreak}

Define the rescaled kernel $\bar{K}_\tau$ as
%
%
\begin{equation}
\label{eq:K-tau-1}
\bar{K}_\tau(s,t) = \sqrt{\tau'(s)\tau'(t)} K_N(\tau(s),\tau(t))
=\sigma_{n,N}K_N(\tau(s),\tau(t)).
\end{equation}
Since $K_N$ and $\bar{K}_\tau$ share the spectrum,
$F_{N,1}(s') = \sqrt{\det(I - \bar{K}_\tau f)}.$

To work out a representation for $\bar{K}_\tau$, apply the $\tau
$-scaling to
$\phi,\psi$ and $S_{N,2}$ to define
%
%
\begin{equation}
\label{eq:phi-psi-tau}
\phi_\tau(s) = \sigma_{n,N}\phi(\tauexpr{s}),\qquad
\psi_\tau(s) = \sigma_{n,N}\psi(\tauexpr{s})
\end{equation}
and
%
%
\begin{equation}
\label{eq:s-tau}
S_\tau(s,t)  = \sigma_{n,N}S_{N,2}(\tauexpr{s},\tauexpr{t})
= (\phi_\tau\diamond\psi_\tau+ \psi_\tau\diamond\phi_\tau)(s,t).
%
\end{equation}
Then we obtain from \eqref{eq:central} that
%
%
\begin{equation}
\label{eq:S-tau-R}
S_\tau^R(s,t) = \sqrt{\tau'(s)\tau'(t)}S_{N,1}(\tau(s),\tau(t))
= S_\tau(s,t) +\psi_\tau(s) (\varepsilon\phi_\tau)(t).
\end{equation}
This, together with \eqref{eq:K-N-kernel} and \eqref{eq:K-tau-1}, leads
to
\[
\bar{K}_\tau(s,t) =
\pmatrix{ I & -\sigma_{n,N}^{-1}\cdot\partial_2\vspace*{3pt}\cr
\sigma_{n,N}\cdot\varepsilon_1 & T
}
S_\tau^R(s,t) + \sigma_{n,N}K^\varepsilon(s,t).
\]
Observe that $\det(I - \bar{K}_\tau f)$ remains unchanged if we divide
the lower left entry by $\sigma_{n,N}$ and multiply the upper right
entry by $\sigma_{n,N}$. Thus, we obtain
%
%
\begin{equation}
\label{eq:F-N-2}
F_{N,1}(s') = \sqrt{\det(I - K_\tau f)}
\end{equation}
with
%
%
\begin{equation}
\label{eq:K-tau-2}
K_\tau(s,t) = (LS_\tau^R)(s,t) + K^\varepsilon(s,t).
\end{equation}

To match the representation \eqref{eq:K-GOE} of $K_{\mathrm{GOE}}$, and to
facilitate later arguments, it is helpful to rewrite $LS_\tau^R$, and
hence $K_\tau$, using
$\tee$. To this end, observe that $\int\psi_\tau= 0$ and~let
%
%
\begin{equation}
\label{eq:beta-N}
\beta_N = \frac{1}{2}\int_{-\infty}^\infty\phi_\tau(s)\,\mathrm{d}s.
\end{equation}
By the identity $(\varepsilon g)(s) = \frac{1}{2}\int g - (\tee
g)(s)$, we
obtain $\varepsilon\phi_\tau= \beta_N - \tee\phi_\tau$ and
$\varepsilon\psi_\tau=
-\tee\psi_\tau$, and so
\[
LS_\tau^R = L(S_\tau-\psi_\tau\otimes\tee\phi_\tau) + \beta_N
L(\psi_\tau\otimes1).
\]
Now $L = \tilde{L} + E$ with $E =
\left({0\atop \varepsilon_1+\tee_1}\enskip{0\atop 0}
\right)
$. Since $2(\varepsilon_1 + \tee_1)$
equals integration over $\mathbb{R}$ in the first argument and
$\int\psi_\tau= 0$, we obtain
\begin{eqnarray*}
LS_\tau^R & =& \tilde{L}(S_\tau- \psi_\tau\otimes\tee\phi_\tau) +
ES_\tau+ \beta_N \tilde{L}(\psi_\tau\otimes1) \\
& =& \tilde{L}(S_\tau- \psi_\tau\otimes\tee\phi_\tau) + \beta_N
L_1(\psi_\tau\otimes1) + \beta_N L_2(1\otimes\psi_\tau).
\end{eqnarray*}
The second equality holds, for $(ES_\tau)_{21} =
\frac{1}{2}\int_{-\infty}^\infty S_\tau(u,t)\,\mathrm{d}t = \beta_N\int
_0^\infty
\psi_\tau(t+z)\,\mathrm{d}z =  \beta_N(\tee\psi_\tau)(t)$. Finally, this gives
$K_\tau$ a similar decomposition to that of $K_{\mathrm{GOE}}$
%
%
\begin{equation}
\label{eq:k-tau-decomp}
K_\tau= \tilde{L}(S_\tau- \psi_\tau\otimes\tee\phi_\tau) +
L_1(\psi_\tau\otimes\beta_N) + L_2(\beta_N\otimes\psi_\tau) +
K^\varepsilon.
\end{equation}

\subsection{Generalized Fredholm determinants}
\label{sec:det}

For any fixed $s_0\in\mathbb{R}$, we are interested in the convergence
rate of $F_{N,1}(s')$ to $F_1(s')$ for all $s'\geq s_0$. In what
follows, we show that this relies on the operator convergence of
$K_\tau$ to $K_{\mathrm{GOE}}$.

First, we note that the determinants in \eqref{eq:F-N-1},
\eqref{eq:F-1-2} and \eqref{eq:F-N-2} are not the usual Fredholm
determinants {}(see, e.g., \cite{lax} for an introduction
to the Fredholm determinant), as the~$\varepsilon$ term on the lower-left
position of the matrix kernels is not of trace class. Tracy and Widom
\cite{tw05} first observed the problem and proposed a solution by
introducing weighted Hilbert spaces and regularized \mbox{2-determinants},
which we adopt here.

Consider the determinant in \eqref{eq:F-N-1}. Let $\tilde\rho$ be a
weight function such that (1) its reciprocal $\tilde\rho^{-1}\in
L^1[0,\infty)$; and (2) $S_{N,1}\in L^2((x',\infty);\tilde\rho)\cap
L^2((x',\infty);\tilde\rho^{-1})$. Then $\varepsilon$:
$L^2((x',\infty);\tilde\rho)\to L^2((x',\infty);\tilde\rho^{-1})$ is
Hilbert--Schmidt and $K_N$ can be regarded as a $2\times2$ matrix
kernel on the space $L^2((x',\infty);\tilde\rho)\oplus
L^2((x',\infty);\tilde\rho^{-1})$. In addition, by the second
condition on~$\tilde\rho$, the diagonal elements of $K_N$ are trace class on
$L^2((x',\infty);\tilde\rho)$ and $L^2((x',\infty);\tilde\rho^{-1}),$
respectively.

For a Hilbert--Schmidt operator $T$ with eigenvalues $\mu_k$, its
regularized 2-determinant~\cite{gohberg00} is defined as $ {\det}_2(I-T)
\equiv\prod_{k}(1-\mu_k)\mathrm{e}^{\mu_k}$. If the diagonal elements of
$T$ are trace class, then we define the generalized Fredholm
determinant for
$T$ as
%
%
\begin{equation}
\label{eq:det}
\det(I-T) = {\det}_2(I-T)\exp(-\trace T).
\end{equation}
As remarked in \cite{tw05}, the definition \eqref{eq:det} is independent
of the choice of $\tilde\rho$ and allows the derivation in \cite{tw98}
that yields \eqref{eq:F-N-1}, \eqref{eq:K-N-kernel} and eventually
\eqref{eq:central}.

Change the domain to $(s',\infty)$ with $s'=\tau^{-1}(x')$ and the
weight function to \mbox{$\rho
= \tilde\rho\circ\tau$}, and abbreviate $L^2((s',\infty);\varrho)$ as
$L^2(\varrho)$ for any suitable $\varrho$. Then, $K_\tau$ and $K_{\mathrm{GOE}}$
are~mem\-bers of the operator class $\cal{A}$ of $2\times2$
Hilbert--Schmidt operator matrices on $L^2(\rho)\oplus L^2(\rho^{-1})$
with trace class diagonal entries. Definition \eqref{eq:det} and
previous derivations in Section \ref{sec:rescale} remain valid.

In order to make the latter argument more explicit, it is convenient to
make a specific choice of the weight function $\rho$. In particular, on
the $s$-scale, we choose
%
%
\begin{equation}
\label{eq:rho}
\rho(s) = 1 + \exp(|s| ). 
\end{equation}
This implies that on the $x$-scale, we specify the
weight function $\tilde\rho= \rho\circ(\tau^{-1})$ as
\[
\tilde\rho(x) = 1 + \exp(| x -
\mu_{n,N}|/{\sigma_{n,N}} ).
\]
It is straightforward to verify that the required conditions are all satisfied.

With rigorous definition of the determinants, we now relate the convergence
of $F_{N,1}$ to $F_1$ to that of $K_\tau$ to $K_{\mathrm{GOE}}$. First of all,
simple manipulation leads to
%
%
\begin{equation}
\label{eq:diff-det}
| F_{N,1}(s') - F_1(s') | \leq\frac{| F_{N,1}^2(s') -
F_1^2(s')|}{F_1(s_0)} = \frac{1}{F_1(s_0)}| \det(I-K_\tau) -
\det(I - K_{\mathrm{GOE}}) | .\vspace*{-2pt}
\end{equation}
To bound the difference between the determinants, we have the following
Lipschitz-type inequality. Here and after, $\|\cdot\|_1$ and $\|\cdot
\|_2$
denote the trace class norm and Hilbert--Schmidt norm, respectively.
\begin{proposition}
\label{prop:det-bd}
Let $A, B \in\mathcal{A}$, and $\det(I-A)$, $\det(I-B)$ defined as
in \eqref{eq:det}. If $\sum_{i=1}^2\|A_{ii} - B_{ii}\|_1 +
\sum_{i\neq j}\|A_{ij} - B_{ij}\|_2\leq1/2$, then
%
%
\begin{equation}
\label{eq:det-bd}
| \det(I-A) - \det(I-B) | \leq
M(B) \Biggl(\sum_{i=1}^2\|A_{ii}-B_{ii}\|_1 + \sum_{i\neq
j}\|A_{ij} -
B_{ij}\|_2 \Biggr) ,\vspace*{-2pt}
\end{equation}
where $M(B) = 2 | \det(I-B)| + 2\exp[2 (1+
\|B\|_2 )^2 + \sum_i\|B_{ii}\|_1 ].$
\end{proposition}
\begin{pf}
\cite{johnstone07}, Proposition 3, established a similar
bound to \eqref{eq:det-bd}, but with $M(B)$ replaced by
\[
C(A,B) = | \mathrm{e}^{-\trace{A}} | \exp\biggl[\frac{1}{2}(1+\|A\|_2 +
\|B\|_2)^2\biggr] + |\operatorname{det}_2(I-B)| \frac{|\mathrm{e}^{-\trace{A}} -
\mathrm{e}^{-\trace{B}}|}{|\trace{A} - \trace{B}|} .\vspace*{-2pt}
\]
We now bound $C(A,B)$ by the above claimed constant $M(B)$.

Observe that for $|x|\leq1/2$, $| \mathrm{e}^x - 1 | \leq
2|x|$. Therefore, when $\sum_{i=1}^2\|A_{ii} - B_{ii}\|_1 + \sum
_{i\neq
j}\|A_{ij} - B_{ij}\|_2\leq1/2$, we have $|\trace{A} - \trace{B}|
\leq\sum_{i=1}^2 \|A_{ii} - B_{ii}\|_1\leq1/2$, which in turn implies
$|\mathrm{e}^{-\trace{A}} - \mathrm{e}^{-\trace{B}}| \leq2|\trace{A} -
\trace{B}||\mathrm{e}^{-\trace{B}}|$. Hence, for the terms in $C(A,B)$, we have
\begin{eqnarray*}
| \mathrm{e}^{-\trace{A}}| & \leq&| \mathrm{e}^{-\trace{B}} -\mathrm{e}^{-\trace{A}}
| + | \mathrm{e}^{-\trace{B}}| \leq| \mathrm{e}^{-\trace{B}}|
(2| \trace{A} - \trace{B} |
+ 1)\\[-2pt]
& \leq&| \mathrm{e}^{-\trace{B}}| \biggl(2\sum_{i}\|A_{ii} -
B_{ii}\|_1 + 1\biggr) \leq2\exp(\|B_{11}\|_1 + \|B_{22}\|_1)\vspace*{-2pt}
\end{eqnarray*}
and
\[
|\operatorname{det}_2(I-B)| \frac{|\mathrm{e}^{-\trace{A}} -
\mathrm{e}^{-\trace{B}}|}{|\trace{A} - \trace{B}|} \leq2 |\operatorname{det}_2(I-B)|
|\mathrm{e}^{-\trace{B}}| = 2|\det(I-B)|.\vspace*{-2pt}
\]
Moreover, we observe that
\begin{eqnarray*}
1 + \|A\|_2 + \|B\|_2 & \leq&1 + 2\|B\|_2 + \|A-B\|_2 \\[-2pt]
& \leq&1 + 2\|B\|_2 + \sum_{i=1}^2\|A_{ii}-B_{ii}\|_1 +
\sum_{i\neq j}\|A_{ij} -
B_{ij}\|_2 \\[-2pt]
&\leq&2 + 2\|B\|_2.\vspace*{-2pt}
\end{eqnarray*}
Plugging all these bounds into $C(A,B)$, we obtain the claimed form
of\vadjust{\goodbreak}
$M(B)$.
\end{pf}

\begin{remark}
Proposition \ref{prop:det-bd} refines \cite{johnstone07}, Proposition
3, by having the leading constant~$M(B)$ of the bound
depend only on $B$, which is important for deriving properties
of the~$C(s_0)$ function later.\vspace*{-2pt}
\end{remark}

\subsection{\texorpdfstring{Decomposition of $K_\tau-K_{\mathrm{GOE}}$}
{Decomposition of K tau - K GOE}}\vspace*{-2pt}
\label{sec:ker-diff}

By Proposition \ref{prop:det-bd}, to prove
Theorem \ref{thm:main} is essentially to control the entrywise
convergence rate of $K_\tau$ to $K_{\mathrm{GOE}}$. To this end, we construct a
telescopic decomposition of $K_\tau- K_{\mathrm{GOE}}$ into sums of simpler
matrix kernels whose entries are more tractable.

To explain the intuition behind the decomposition, we introduce
constants $\tmu_{n,N}$ and~$\tsigma_{n,N}$ as
%
%
\begin{equation}
\label{eq:lag-center-scale}
\tmu_{n,N} =\bigl (\sqrt{n_+} + \sqrt{N_+} \bigr)^2,\qquad
\tsigma_{n,N} = \bigl(\sqrt{n_+} +
\sqrt{N_+} \bigr)\biggl ( \frac{1}{\sqrt{n_+}} + \frac{1}{\sqrt{N_+}}
\biggr)^{1/3}.\vspace*{-2pt}
\end{equation}
In \cite{nek06}, it was shown that $(\mu_{n, N}, \sigma_{n,N}) =
(\tmu_{n-1,N-1}, \tsigma_{n-1,N-1})$ is `optimal' for $\psi_\tau$
in the
sense that $| \psi_\tau- G | = \Oh{N^{-2/3}}$, but suboptimal
for $\phi_\tau$ as $| \phi_\tau- G| = \Oh{N^{-1/3}}$. However,
later in Proposition \ref{prop:laguerre}, we will show that $|
\phi_\tau- G - \Delta_N G' | = \Oh{N^{-2/3}}$ for
%
%
\begin{equation}
\label{eq:Delta-N}
\Delta_N = \frac{\tmu_{n-1,N-1}-\tmu_{n-2,N}}{\tsigma_{n-2,N}} =
\Oh{N^{-1/3}}.\vspace*{-2pt}
\end{equation}
(For a proof, see Section \ref{subsec:a-constants}.)
These bounds suggest that, in the decomposition, we align~$\psi_\tau$
with $G$, and $\phi_\tau$ with $G + \Delta_NG'$.

Let $G_N = G + \Delta_N G'$ and $S_{A_N} = G\diamond G_N +
G_N\diamond G$. We obtain
\[
S_{A_N} - G\otimes\tee G_N = S_A - G\otimes\tee G\vspace*{-2pt}
\]
for
\[
\int_0^\infty G(s+z)G'(t+z)+G'(s+z)G(t+z)\,\mathrm{d}z = \int_0^\infty
\frac{\mathrm{d}}{\mathrm{d}z}[G(s+z)G(t+z)]\,\mathrm{d}z = -G(s)G(t).\vspace*{-2pt}
\]
This, together with \eqref{eq:K-GOE} and \eqref{eq:k-tau-decomp},
leads to the decomposition
%
%
\begin{eqnarray}\label{eq:kern-diff}
K_\tau-K_{\mathrm{GOE}} & = &\tilde{L}(S_\tau- S_{A_N}) + \tilde{L}(G\otimes
\tee G_N - \psi_\tau\otimes\tee\phi_\tau) \nonumber
\\[-9pt]
\\[-9pt]
&&{}  + L_1\bigl(\psi_\tau\otimes\beta_N - G\otimes
\tfrac{1}{\sqrt{2}}\bigr) + L_2\bigl(\beta_N\otimes\psi_\tau-
\tfrac{1}{\sqrt{2}}\otimes G\bigr).
\nonumber\vspace*{-2pt}
\end{eqnarray}

\subsection{Laguerre asymptotics and operator bounds}\vspace*{-2pt}
\label{sec:lemma}

Here we collect a set of intermediate results to be used repeatedly in
the proof of Theorem~\ref{thm:main}.

To start with, we consider the asymptotics of $\phi_\tau$ and $\psi
_\tau$
and their derivatives. Recalling that $G = \frac{1}{\sqrt{2}}\Ai$
and $G_N
= G + \Delta_N G'$, we have the following.
\begin{proposition}
\label{prop:laguerre}
Let $\phi_\tau$, $\psi_\tau$ and $\Delta_N$ be defined as in
\eqref{eq:phi-psi-tau} and \eqref{eq:Delta-N}. Assume that
\eqref{eq:ch-var} holds, and that as $N\to\infty$, $n =
n(N)\to\infty$ with $n/N\to\gamma\in[1,\infty)$. Then,\vadjust{\goodbreak} for any given~$s_0$, there exists an integer $N_0(s_0,\gamma)$ such that when
$N\geq N_0(s_0,\gamma)$, for all $s\geq s_0$,
%
\begin{eqnarray}\label{eq:psi-psi'-bd}
| \psi_\tau(s)| ,| \psi_\tau'(s) | & \leq&
C(s_0)\exp(-s), \\\label{eq:phi-phi'-bd}
| \phi_\tau(s)| ,| \phi_\tau'(s) | & \leq&
C(s_0)\exp(-s), \\\label{eq:psi'-bd}
| \psi_\tau(s) - G(s) | ,| \psi_\tau'(s) - G'(s) | &
\leq& C(s_0)N^{-2/3}\exp(-s), \\\label{eq:phi'-bd}
| \phi_\tau(s) - G_N(s) | ,| \phi_\tau'(s) -
G_N'(s) | & \leq&
C(s_0)N^{-2/3}\exp(-s),
\end{eqnarray}
where $C(\cdot)$ is continuous and non-increasing.
\end{proposition}

Integrating these bounds over $[s,\infty)$, we know that they remain
valid if we replace~$\psi_\tau, \phi_\tau, G$ and $G_N$ with $\tee
\psi_\tau, \tee\phi_\tau, \tee G$ and $\tee G_N$ on the left-hand sides.
The proof of Proposition \ref{prop:laguerre} involves careful
Liouville--Green analysis on the solution of certain differential
equations and will be discussed in detail later in Section
\ref{sec:lg}.

On the other hand, for $G$ and $G_N$, we have the following bounds from
\cite{olver74}, page~394. Note that the bounds for $G_N$ and $G_N'$ do not
depend on $N$, for $\Delta_N$ is uniformly bounded.%

\begin{lemma}
\label{lemma:airy}
Fix $\beta> 0$ and $k\geq0$. Then, for all $s\geq s_0$,
\[
|s^k G(s)|,|s^k G_N(s)|,|s^k G'(s)|,|s^k G'_N(s)| \leq C(s_0) \exp
(-\beta s),
\]
where $C(s_0)$ is continuous and non-increasing.
\end{lemma}

For a proof of the lemma, see \cite{twacclong}. Integrating the bounds
for $|G|$ and $|G_N|$ over $[s,\infty)$, we obtain that $|\tee G|$ and
$|\tee G_N|$ are also bounded by $C(s_0)\mathrm{e}^{-\beta s}$.


For a later operator convergence argument, we will need simple bounds
for certain
norms of operator $D\dvtx L^2(\rho_1)\to L^2(\rho_2)$ with kernel $D(u,v) =
\aa(u)\beta(v)(a\diamond b)(u,v)$, where
$\{\rho_1,\rho_2\}\subset\{\rho,\rho^{-1}\}$ with $\rho$ given in
\eqref{eq:rho}. In particular, we have
\begin{lemma}[(\cite{johnstone07})]
\label{lemma:hs-tr-bd}
Let $D\dvtx L^2(\rho_1) \to L^2(\rho_2)$ have kernel $D(u,v) =
\aa(u)\beta(v)(a\diamond b)(u,v)$. Suppose that
$\{\rho_1,\rho_2\}\subset\{\rho,\rho^{-1}\}$ and that, for $u\geq s'$,
%
%
\begin{equation}
\label{eq:hs-tr-bd1}
|\aa(u)|\leq\aa_0\mathrm{e}^{\aa_1 u}, \qquad|\beta(u)|\leq\beta
_0\mathrm{e}^{\beta_1
u}, \qquad|a(u)|\leq a_0\mathrm{e}^{-a_1 u}, \qquad|b(u)|\leq b_0\mathrm{e}^{-b_1 u},
\end{equation}
with $a_1-\aa_1, b_1-\beta_1 \geq1$. Then the Hilbert--Schmidt norm
satisfies
%
%
\begin{equation}
\label{eq:D-bd-1}
\|D\|_2 \leq C \frac{\aa_0 \beta_0 a_0
b_0}{a_1+b_1} \exp[-(a_1+b_1-\aa_1-\beta_1)s' + |s'|
],
\end{equation}
where $C = C(a_1,\aa_1, b_1, \beta_1)$. If $\rho_1 = \rho_2$, the trace
norm $\|D\|_1$ satisfies the same bound.
\end{lemma}

\subsection{\texorpdfstring{Operator convergence: Proof of Theorem \protect\ref{thm:main}}
{Operator convergence: Proof of Theorem 1}}
\label{sec:conv}

Abbreviate the terms in the decomposition \eqref{eq:kern-diff} as
\[
K_\tau- K_{\mathrm{GOE}} = \delta^R + \delta_0^F + \delta_1^F +
\delta_2^F.
\]
We work out below entrywise bounds for each of these $\delta$ terms and
then apply Proposition~\ref{prop:det-bd} to complete the proof of Theorem
\ref{thm:main}. In what follows, we use the abbreviation~$D^{(k)}f$, $k
= -1, 0, 1$ to denote $\tee f$, $f$ and $f'$, respectively. Moreover, the
unspecified norm $\|\cdot\|$ denotes the Hilbert--Schmidt norm $\|
\cdot\|_2$
for off-diagonal entries and trace class norm $\|\cdot\|_1$ for diagonal
ones.

\subsubsection*{$\delta^R$ term} Recall that $\delta^R = \tilde
{L}(S_\tau-
S_{A_N})$ with $S_\tau= \phi_\tau\diamond\psi_\tau+ \psi_\tau
\diamond
\phi_\tau$ and $S_{A_N} = G_N\diamond G + G \diamond G_N$. Regardless
of the signs, we have the following unified expression for the entries
of $\delta^R$:
%
%
\begin{eqnarray}
\label{eq:delta-R-entries}
(\delta^R)_{ij}& = & D^{(k)}(\phi_\tau- G_N)\diamond
D^{(l)}\psi_\tau+ D^{(k)}G_N\diamond D^{(l)}(\psi_\tau- G)\nonumber
\\[-8pt]
\\[-8pt]
&&{} + D^{(k)}(\psi_\tau- G)\diamond D^{(l)}\phi_\tau+
D^{(k)}G\diamond D^{(l)}(\phi_\tau- G_N),
\nonumber
\end{eqnarray}
for $i,j\in\{1, 2\}$, $k\in\{-1, 0\}$ and $l\in\{0, 1\}$. By
Proposition \ref{prop:laguerre} and Lemma \ref{lemma:airy}, we find
that for any of the four terms in \eqref{eq:delta-R-entries},
condition \eqref{eq:hs-tr-bd1} is satisfied with $\aa_0 = \beta_0 =
1$, $\aa_1 = \beta_1 = 0$, $a_1 = b_1 = 1$ and $\{a_0, b_0\} = \{C(s_0),
C(s_0)N^{-2/3}\}$.
So Lemma \ref{lemma:hs-tr-bd} implies
%
%
\begin{equation}
\label{eq:delta-R-bd}
\|(\delta^R)_{ij}\| \leq C(s_0)N^{-2/3}\exp( -2s' + |s'|
).
\end{equation}
By a simple triangle inequality, we can choose $C(s_0)$ in the last
display as the sum of products of continuous and non-increasing
functions, which can be seen from the term
$(\aa_0\beta_0a_0b_0)/(a_1+b_1)$ in \eqref{eq:D-bd-1}. Moreover, the
term $C$ in \eqref{eq:D-bd-1} is a universal constant for fixed $a_1,
\aa_1, b_1$ and $\beta_1$ here. Hence, the final $C(s_0)$ function
remains continuous and non-increasing.

\subsubsection*{Finite rank terms} For a rank one operator $a\otimes
b\dvtx
L^2(\rho_1)\to L^2(\rho_2)$ with kernel $a(s)b(t)$, its norm is
\[
\| a\otimes b \| = \|a\|_{2,\rho_2}\|b\|_{2,\rho_1^{-1}}.
\]
Here, the norm can be either trace class or Hilbert--Schmidt, since the two
agree for rank one operators. In addition, for any $\varrho$,
$\|a\|^2_{2,\varrho} = \int_{s'}^\infty|a(s)|^2 \varrho(s)\,\mathrm{d}s$. Now
consider matrices of rank one operators on $L^2(\rho)\otimes
L^2(\rho^{-1})$. Write $\|\cdot\|_+$ and $\|\cdot\|_-$ for
$\|\cdot\|_{2,\rho}$ and $\|\cdot\|_{2,\rho^{-1}}$,
respectively. \cite{johnstone07}, equation (213) gives the following bound
%
%
\begin{equation}
\label{eq:rankone-bd}
\pmatrix{
\|a_{11}\otimes b_{11} \|_1 & \|a_{12}\otimes b_{12} \|_2 \cr
\|a_{21}\otimes b_{21} \|_2 & \|a_{22}\otimes b_{22} \|_1
}
\leq
\pmatrix{
\|a_{11}\|_+\|b_{11}\|_- & \|a_{12}\|_+\|b_{12}\|_+ \cr
\|a_{21}\|_-\|b_{21}\|_- & \|a_{22}\|_-\|b_{22}\|_+
}
.
\end{equation}

First consider $\delta^F_0$. We reorganize it as
\begin{eqnarray*}
\delta^F_0 &=& -\tilde{L}(\psi_\tau\otimes\tee\phi_\tau- G\otimes
\tee G_N)\\
 &=& -\tilde{L}[\psi_\tau\otimes\tee(\phi_\tau- G_N) +
(\psi_\tau- G)\otimes\tee G_N] = \delta^{F,1}_0 + \delta^{F,2}_0.
\end{eqnarray*}
The entries of $\delta^{F,i}_0$, $i=1, 2$, are all of the form
$a\otimes b$,
with $a$ and $b$ chosen from $D^{(k)}\psi_\tau$, $D^{(k)}(\phi_\tau
- G_N)$, $D^{(k)}(\psi_\tau- G)$ and $D^{(k)}G_N$, for $k\in\{-1, 0,
1\}$.\vadjust{\goodbreak}

Observe that for $\eta\geq2$ we have
%
%
\begin{equation}
\label{eq:eta}
\int_{s'}^\infty\exp(-\eta s)\rho^{\pm1}(s)\,\mathrm{d}s \leq\frac{4}{\eta-
1}\exp(-\eta s' \pm|s'|)\leq\frac{8}{\eta}\exp(-\eta s'+|s'|).
\end{equation}
Together with Proposition \ref{prop:laguerre} and Lemma
\ref{lemma:airy}, this implies
\begin{eqnarray*}
\bigl\| D^{(k)}\psi_\tau\bigr\|^2_\pm, \bigl\|D^{(k)} G_N\bigr\|^2_\pm& \leq&
C(s_0)\exp(-2s'+|s'|),\\
\bigl\|D^{(k)}(\psi_\tau- G)\bigr\|^2_\pm, \bigl\|D^{(k)}(\phi_\tau- G_N)\bigr\|^2_\pm
& \leq&
C(s_0)N^{-4/3}\exp(-2s'+|s'|).
\end{eqnarray*}
These bounds, together with the triangle inequality and \eqref
{eq:rankone-bd}, yield
\begin{eqnarray*}
\|(\delta_0^F)_{11}\|_1
& \leq&\| \psi_\tau\otimes\tee(\phi_\tau-G_N) \|_1
+ \| (\psi_\tau-G)\otimes\tee G_N \|_1 \\
& \leq&\|\psi_\tau\|_+ \|\tee(\phi_\tau-G_N)\|_-
+ \|\psi_\tau- G\|_+ \| \tee G_N \|_- \\
& \leq& C(s_0)N^{-2/3}\exp(-2s'+|s'|).
\end{eqnarray*}
Similarly, we obtain the bounds for the other entries. In summary, we have
%
%
\begin{equation}
\label{eq:delta-0-bd}
\|(\delta^F_0)_{ij}\| \leq C(s_0)N^{-2/3}\exp(-2s' +
|s'| ).
\end{equation}

Switch to $\delta^F_1$ and $\delta^F_2$. Recall that $\delta^F_1 =
L_1(\psi_\tau\otimes\beta_N - G \otimes\frac{1}{\sqrt{2}})$ and
$\delta^F_2 = L_2(\beta_N\otimes\psi_\tau- \frac{1}{\sqrt
{2}}\otimes
G)$. Due to their similarity, we take $\delta^F_1$ as an example and the
same analysis applies to $\delta^F_2$ with obvious modification. We
further decompose
$\delta^F_1$as
\[
\delta^F_1 = L_1\bigl[ (\psi_\tau- G)\otimes\beta_N + G\otimes
\bigl(\beta_N - \tfrac{1}{\sqrt{2}} \bigr) \bigr]. 
\]
By \eqref{eq:rankone-bd}, the essential elements we need to bound are
$\|D^{(k)}(\psi_\tau- G)\|_\pm$, $\|D^{(k)}G\|_\pm$ and $\|1\|_-$ for
$k = -1$ and $0$. The bounds related to $D^{(k)}(\psi_\tau- G)$ have
already been obtained. For the other two terms, \eqref{eq:eta} and Lemma
\ref{lemma:airy} give
\[
\bigl\|D^{(k)}G\bigr\|^2_\pm\leq C(s_0)\exp(-2s' + |s'| )
\]
and
\[
\|1\|_-^2 = \int_{s'}^\infty[1 +
\exp(|s|) ]^{-1}\,\mathrm{d}s \leq\int_{-\infty}^\infty\exp(-|s|)\,\mathrm{d}s
\leq2.
\]
Since $\beta_N - \frac{1}{\sqrt{2}} = \Oh{N^{-1}}$ (for a proof, see
Section
\ref{subsec:a-constants}), we have
\begin{eqnarray*}
\|(\delta^F_1)_{11}\|_1 &\leq&\|(\psi_\tau- G)\otimes\beta_N\|_1
+\bigl\| G\otimes\bigl(\beta_N - 1/\sqrt{2}\bigr)\bigr\|_1 \\
&\leq& \|(\psi_\tau- G)\|_+ \|\beta_N\|_- + \|G\|_+ \bigl\|\beta_N -
1/\sqrt{2}\bigr\|_- \\
&\leq& C(s_0)N^{-2/3}\exp(-s' + |s'|/2 ) +
C(s_0)N^{-1}\exp(-s' + |s'|/2 ) \\
&\leq& C(s_0)N^{-2/3}\exp(-s'/2).
\end{eqnarray*}
In a similar vein, the same bound can be obtained for
$\|(\delta^{F}_1)_{12}\|_2$ and entries of $\delta^F_2$. Therefore, we
conclude that
%
%
\begin{equation}
\label{eq:delta-1-2-bd}
\|(\delta^F_1)_{ij}\|, \|(\delta^F_2)_{ij}\| \leq
C(s_0)N^{-2/3}\exp(-s'/2 ).
\end{equation}

Now we prove Theorem \ref{thm:main}.

\begin{pf*}{Proof of Theorem \ref{thm:main}}
By the decomposition \eqref{eq:kern-diff} and bounds
\eqref{eq:delta-R-bd}, \eqref{eq:delta-0-bd} and
\eqref{eq:delta-1-2-bd}, the triangle inequality gives the following
bound for the norm of each entry in $K_\tau-K_{\mathrm{GOE}}$:
\[
\|(K_{\tau} -K_{\mathrm{GOE}})_{ij} \|\leq C(s_0)N^{-2/3}\exp(-s'/2).
\]
We then apply Proposition \ref{prop:det-bd} with $A = K_\tau$ and $B =
K_{\mathrm{GOE}}$ to get
%
%
\begin{equation}\label{eq:det-diff-bd}
| \det(I-K_\tau) - \det(I-K_{\mathrm{GOE}}) | \leq
M(K_{\mathrm{GOE}})C(s_0)N^{-2/3}\exp(-s'/2),
\end{equation}
where $M(K_{\mathrm{GOE}}) = 2\det(I-K_{\mathrm{GOE}}) + 2\exp\{2 ( 1 +
\|K_{\mathrm{GOE}}\|_2 )^2 + \sum_i \|K_{\mathrm{GOE},ii}\|_1\}$.

For the first term in $M(K_{\mathrm{GOE}})$, we have $\det(I-K_{\mathrm{GOE}}) =
F_1^2(s') \leq1$. On the other hand, we have
\[
\|K_{\mathrm{GOE}}\|_2 \leq\sum_{i,j}\|(K_{\mathrm{GOE}})_{ij}\|_2 \leq
\sum_i\|(K_{\mathrm{GOE}})_{ii}\|_1 + \sum_{i\neq j}\|(K_{\mathrm{GOE}})_{ij}\|_2.
\]
In principle, one can show that, for each $(i,j)$, $\|(K_{\mathrm{GOE}})_{ij}\|
\leq C(s_0)$, with $C(s_0)$ continuous and non-increasing. Take
$\|(K_{\mathrm{GOE}})_{11}\|_1$ as an example. Let $H_\tau$ and $G_\tau$ be
Hilbert--Schmidt operators with kernels $\phi_\tau(x+y)$ and
$\psi_\tau(x+y),$ respectively, then as an operator
\[
(K_{\mathrm{GOE}})_{11} = H_\tau G_\tau+ G_\tau H_\tau+ G\otimes
\tfrac{1}{\sqrt{2}} - G\otimes\tee G.
\]
Since $\|AB\|_1 \leq\|A\|_2\|B\|_2$, we have
\[
\|(K_{\mathrm{GOE}})_{11}\|_1 \leq2 \|H_\tau\|_2\|G_\tau\|_2 +
\tfrac{1}{\sqrt{2}}\|G\|_{2,\rho} \|1 \|_{2,\rho^{-1}} +
\|G\|_{2,\rho}\|\tee G\|_{2,\rho^{-1}}.
\]
Each norm on the right-hand side of the last inequality is the square root
of an
integral of a positive function on $(s',\infty)$ or $(s',\infty)^2$ that
is bounded by the corresponding integral over $(s_0,\infty)$ or
$(s_0,\infty)^2$, which in turn is continuous and non-increasing in
$s_0$. Hence, $\|(K_{\mathrm{GOE}})_{11}\|_1 \leq C(s_0)$. A similar argument
applies to other entries. So, we can control $M(K_{\mathrm{GOE}})$ by a
continuous and non-increasing $C(s_0)$. Finally, we complete the proof
by noting \eqref{eq:diff-det} and the fact that $1/F_1(s_0)$ is
continuous and non-increasing.
\end{pf*}

\section{The smallest eigenvalue}
\label{sec:proof-thm2}

This section is dedicated to the proof of Theorem \ref{thm:small}.

Recall that two key components in the proof of Theorem \ref{thm:main} were:
(1) determinantal representations for both the finite and the limiting
distributions; (2) a closed-form formula for the finite sample kernel
that yields a convenient decomposition of its difference from the
limiting kernel.

In what follows, we first establish the rate of convergence for matrices
with even dimensions. This is achieved by working out the above two
components in the case of the smallest eigenvalue. Then, we prove weak
convergence for matrices with odd dimensions using an interlacing
property of the singular values.

\subsection{Determinantal formula}
\label{sec:det-small}

As before, we follow RMT notation to replace $p$ with $N$, and identify
LOE($N,\aa$) with eigenvalues of $A\sim W_N(I,n)$ by \eqref{eq:ch-var}.

Assume that $N$ is even. For the smallest eigenvalue $x_N$, for any
$x'\geq0$, \cite{tw98} gives
%
%
\begin{equation}
\label{eq:G-N-N}
1 - \tilde{F}_{N,N}(x') = P\{x_N > x'\} = \sqrt{\det(I - K_N\chi)},
\end{equation}
where $\chi= \id_{0\leq x\leq x'}$ and $K_N$ is given in
\eqref{eq:K-N-kernel}.

Due to a nonlinear transformation to be introduced, the formula
\eqref{eq:widom-S1} that we previously used to represent $S_{N,1}$, the
key component in $K_N$, is not most appropriate here. Instead, we find
an alternative (yet equivalent) formula given in \cite{adler00}, Proposition
4.2, more convenient. Indeed, let
%
%
\begin{equation}\label{eq:phi-bar-k}
\bar{\phi}_k(x;\aa) = (-1)^k\sqrt{\frac{a_N}{2}}\phi_k(x;\aa)
x^{-1/2}\id_{x\geq0},
\end{equation}
with $a_N = \sqrt{N(N+\aa)}$. Then \cite{adler00}, Proposition 4.2,
asserts that
%
%
\begin{equation}\label{eq:S-N-1-alt}
S_{N,1}(x,y;\aa) = \sqrt{\frac{y}{x}}S_{N-1,2}(x,y;\aa+ 1) +
\sqrt{\frac{N-1}{N}}\bar{\phi}_{N-1}(x;\aa+1)
(\varepsilon\bar{\phi}_{N-2})(y;\aa+1).
\end{equation}
We write out the explicit dependence of these kernels on the parameter
$\aa$ as they are different on the two sides of the equation. As a
comparison, the previous representation~\eqref{eq:central} could be
rewritten as
\[
S_{N,1}(x,y;\aa) = S_{N,2}(x,y;\aa) +
\bar{\phi}_{N-1}(x;\aa+ 1)(\varepsilon\bar{\phi}_N)(y;\aa- 1).
\]
Its equivalence to \eqref{eq:S-N-1-alt} is given in the Appendix of
\cite{adler00}.

Now, introduce the nonlinear transformation
%
%
\begin{equation}
\label{eq:trans-small}
\pi(s) = \exp(\nu_{n,N}^- - s\tau_{n,N}^-),
\end{equation}
where $\nu_{n,N}^-$ and $\tau_{n,N}^-$ are the rescaling constants in
\eqref{eq:small-center-scale}, with $p$ replaced by $N$. Incorporating
the transformation into $K_N$, we define
%
%
\begin{equation}
\label{eq:K-pi-kern}
\bar{K}_\pi(s,t) = \sqrt{\pi'(s)\pi'(t)} K_N(\pi(s), \pi(t)).
\end{equation}
Let $F_{N,N}$ be the distribution of $(\log{x_N} -
\nu_{n,N}^-)/\tau_{n,N}^-$. Fix $s_0$, for any $s' = \pi^{-1}(x')
\geq
s_0$ and $f = \id_{s\geq s'}$, since\vadjust{\goodbreak} $\det(I-K_N\chi) =
\det(I-\bar{K}_\pi f)$, we obtain $1 - F_{N,N}(-s') =
\sqrt{\det(I-\bar{K}_\pi f)}$. Thinking of $K_\pi$ as a Hilbert--Schmidt
operator with trace class diagonal entries on
$L^2([s',\infty);\rho)\oplus L^2([s',\infty);\rho^{-1})$ for proper
weight function $\rho$, we can drop $f$.

Now consider the
representation of $\bar{K}_\pi$. For $b_N = \sqrt{(N-1)/N}$, let
%
%
\begin{equation}
\label{eq:phi-psi-pi}
\phi_\pi(s) = -\sqrt{b_N}\pi'(s)\bar{\phi}_{N-2}\bigl(\pi(s);\aa+ 1\bigr),
\qquad
\psi_\pi(s) = \sqrt{b_N}\pi'(s)\bar{\phi}_{N-1}\bigl(\pi(s);\aa+1\bigr).
\end{equation}
Using \cite{pjfbook}, Proposition 5.4.2, we obtain
\[
S_{N-1,2}\bigl(\pi(s),\pi(t);\aa+1\bigr) = (\pi'(s)\pi'(t))^{-1/2}
(\phi_\pi\diamond\psi_\pi+ \psi_\pi\diamond\phi_\pi)(s,t).
\]
On the other hand, simple manipulation yields that the second term in
\eqref{eq:S-N-1-alt}, with $x = \pi(s)$ and $y = \pi(t)$, equals
$(-\pi'(s))^{-1}\psi_\pi(s)(\varepsilon\phi_\pi)(t)$. Thus,
$S_{N,1}(\pi(s),\pi(t)) = (-\pi'(s))^{-1} S_\pi^R(s,t)$
with
%
%
\begin{equation}
\label{eq:S-pi-R}
S_\pi^R(s,t) = (\phi_\pi\diamond\psi_\pi+ \psi_\pi\diamond
\phi_\pi)(s,t) + (\psi_\pi\otimes\varepsilon\phi_\pi)(s,t).
\end{equation}
In addition, we have
\begin{eqnarray*}
(-\partial_2S_{N,1})(\pi(s),\pi(t)) & =& \frac{-\partial_t
S_{N,1}(\pi(s),\pi(t))}{\partial_t \pi(t)}
= \frac{-1}{\pi'(s)\pi'(t)}\cdot[-\partial_2 S_\pi^R(s,t)],\\
(\varepsilon_1 S_{N,1})(\pi(s),\pi(t)) & =& \int_0^\infty
\varepsilon(\pi(s)-z)S_{N,1}(z,\pi(t))\,\mathrm{d}z \\
& =& \int_{-\infty}^\infty\varepsilon(s-u)S_{N,1}(\pi(u), \pi
(t))\pi'(u)\,\mathrm{d}u =
-(\varepsilon_1 S_\pi^R)(s,t).
\end{eqnarray*}
Supplying these equations to \eqref{eq:K-N-kernel}, we obtain that
\[
\bar{K}_\pi(\pi(s),\pi(t)) = U(s)(LS_\pi^R + K^\varepsilon)(s,t)U^{-1}(t)
\]
with $U(s) = \operatorname{diag}(1/\sqrt{-\pi'(s)},-\sqrt{-\pi
'(s)})$. Observe
that $\operatorname{det}(I-\bar{K}_\pi)$ remains unchanged if we premultiply
$\bar{K}_\pi$ with $U^{-1}(s_0)$ and postmultiply it with
$U(s_0)$. Denoting the resulting kernel by $K_\pi$, we obtain that
%
%
\begin{equation}
\label{eq:k-pi-entries-log}
K_\pi(s,t) = Q_N(s)(LS_\pi^R + K^\varepsilon)(s,t)Q_N^{-1}(t)
\end{equation}
with $Q_N(s) = U^{-1}(s_0)U(s) =
\operatorname{diag}(\sqrt{\pi'(s_0)/\pi'(s)},\sqrt{\pi'(s)/\pi
'(s_0)})$ and that
$1 -F_{N,N}(-s') = \sqrt{\operatorname{det}(I - K_\pi)}$.

Recall that $G_1(-s') = 1 - F_1(s')$. So, $F_{N,N}(-s') - G_1(-s') =
F_1(-s') - [1 - F_{N,N}(-s')]$. Similar to \eqref{eq:diff-det}, we
obtain
\[
|F_{N,N}(-s') - G_1(-s')| \leq\frac{1}{F_1(s_0)} |\det(I-K_\pi) -
\det(I-K_{\mathrm{GOE}})|.
\]
Thus, as in the case of the largest eigenvalue, by Proposition
\ref{prop:det-bd}, to prove Theorem \ref{thm:small} is to control the
entrywise norm of $K_\pi- K_{\mathrm{GOE}}$. For this purpose, a convenient
decomposition of $K_\pi- K_{\mathrm{GOE}}$ is crucial, to which we now turn.

\subsection{Kernel difference decomposition}
\label{sec:ker-diff-small}

We derive below a decomposition of $K_\pi- K_{\mathrm{GOE}}$. Despite the
differences in actual formulas, the general guideline of the decomposition
is the same as that in Section \ref{sec:ker-diff}.

To start with, we rewrite \eqref{eq:k-pi-entries-log} using the right
tail integration operator $\tee$. To this end, observe that $\int
\psi_\pi= 0$ and that
\[
\tilde\beta_N\,{=}\,\frac{1}{2}\int_{-\infty}^\infty\!\phi_\pi(s)\,\mathrm{d}s\,{=}\,
\frac{(N\,{-}\,1)^{1/4}(n\,{-}\,1)^{1/4}}{2^{(n-N)/2}(N\,{-}\,1)}
\frac{\Gamma( ({N\,{+}\,1})/{2} )}{\Gamma(
{n}/{2} )}
\biggl[\frac{\Gamma(n\,{-}\,1)}{\Gamma(N\,{-}\,1)} \biggr]^{1/2}\!
\,{=}\, \frac{1}{\sqrt{2}}\,{+}\,\Oh{N^{-1}}.
\]
By the same argument that leads to \eqref{eq:k-tau-decomp}, we obtain
\[
K_\pi(s,t) = Q_N(s) (K^R_\pi+ K^F_{\pi,1} + K^F_{\pi,2} +
K^\varepsilon) (s,t)
Q_N^{-1}(t),
\]
with the unspecified components given by
\[
K^R_{\pi} = \tilde{L}(S_\pi- \psi_\pi\otimes
\tee\phi_\pi),
\qquad K^F_{\pi,1} = L_1(\psi_\pi\otimes\tilde\beta_N),
\qquad K^F_{\pi,2} = L_2(\tilde\beta_N\otimes\psi_\pi).
\]

Define $\tilde\Delta_N = (\nu_{n,N}^- -\nu_{n-1,N-1}^-)/\tau_{n-1,N-1}^-
= \Oh{N^{-1/3}}$ and $\tilde{G}_N = G + \tilde\Delta_N G'$. For
$\tilde{S}_{A_N} = G\diamond\tilde{G}_N + \tilde{G}_N\diamond G$, we
have $\tilde{S}_{A_N} - G\otimes\tee\tilde{G}_N = S_A - G\otimes
\tee
G$. Abbreviate the terms in \eqref{eq:K-GOE} as
\[
K_{\mathrm{GOE}} = K^R + K_1^F + K_2^F + K^\varepsilon.
\]
Then,
\begin{eqnarray*}
K^R_\pi- K^R & = &\tilde{L}(S_\pi- S_A - \psi_\pi\otimes\tee\phi
_\pi
+ G\otimes\tee G) \\
& =& \tilde{L}(S_\pi- \tilde{S}_{A_N} ) -
\tilde{L}(\psi_\pi\otimes\tee\phi_\pi- G\otimes\tee\tilde G_N)
= \delta^{R,I} + \delta^F_0.
\end{eqnarray*}
Further define
\begin{eqnarray*}
\delta^{R,D}(s,t) & =& Q_N(s)K^R_\pi(s,t) Q_N^{-1}(t) - K^R_\pi
(s,t),\\
\delta^F_i(s,t) &=& Q_N(s)K^F_{\pi,i}(s,t)Q_N^{-1}(t) -
K^F_i(s,t),\qquad i = 1,2,\\
\delta^\varepsilon(s,t) &=& Q_N(s)K^\varepsilon(s,t) Q_N^{-1}(t) -
K^\varepsilon(s,t).
\end{eqnarray*}
Our final decomposition of $K_\pi- K_{\mathrm{GOE}}$ is
%
%
\begin{equation}
\label{eq:small-ker-decomp}
K_\pi- K_{\mathrm{GOE}} = \delta^{R,D} + \delta^{R,I} + \delta^F_0 +
\delta^F_1 + \delta^F_2 + \delta^\varepsilon.
\end{equation}


We remark that Proposition \ref{prop:laguerre} remains valid if we
replace $\phi_\tau$ and $\psi_\tau$ with $\phi_\pi$ and $\psi
_\pi$,
respectively. The proof is similar to that to be
presented in Section \ref{sec:lg} for Proposition
\ref{prop:laguerre}. 
With these estimates, for each term in \eqref{eq:small-ker-decomp}, we
apply Lemma \ref{lemma:hs-tr-bd} to bound their entrywise norms as in
Section~\ref{sec:conv}. This completes the proof of the rate of
convergence part in Theorem \ref{thm:small}.

\subsection{Weak convergence in the odd $N$ case}

We now establish weak convergence to the reflected Tracy--Widom
law in the odd $N$ case. This is achieved by employing an interlacing
property of the singular values. The strategy follows from
\cite{sosh02}, Remark 5.

Assume that $N$ is odd and $n - 1\geq N$. Let ${X}_{N+1}$ be an
$(n+1)\times(N+1)$ matrix with i.i.d. $N(0,1)$ entries and ${X}_N$ the
$n\times N$ matrix obtained by deleting the last row and the last column
of $X_{N+1}$. Denote the smallest singular values of ${X}_{N+1}$ and
$X_N$ by $\iota_{N+1}$ and $\iota_{N}$, respectively. We apply
\cite{hj85}, Theorem 7.3.9, twice to obtain that $\iota_{N} \leq
\iota_{N+1}$. Repeat the deletion operation on $X_N$ to obtain the
$(n-1)\times(N-1)$ matrix $X_{N-1}$ and denote its smallest singular
value by $\iota_{N-1}$. Then we obtain the `sandwich' relation:
$\iota_{N-1} \leq\iota_N \leq\iota_{N+1}$.

Observe that for $k = N-1, N$ and $N+1$, $X_k'X_k$ are white Wishart
matrices with the smallest eigenvalues $x_k = \iota_k^2$. In addition,
as $N\to\infty$ and $n/N\to\gamma> 1$,
\[
(\nu_{n,N}^- - \nu_{n-1,N-1}^-)/{\tau_{n-1,N-1}^-} =
\Oh{N^{-1/3}}\quad\mbox{and} \quad
{\tau_{n,N}^-}/{\tau_{n-1,N-1}^-} = 1 + \Oh{N^{-1}}.
\]
They together imply that the weak limits for the odd $N$ and the even
$N$ sequences must be the same. This completes the proof of Theorem
\ref{thm:small}.

\section{Laguerre polynomial asymptotics}
\label{sec:lg}

In this section, we complete the proof of Proposition
\ref{prop:laguerre}. The proof has the following components. First, we
take the Liouville--Green approach to analyze an intermediate function
that is connected to both $\phi_\tau$ and $\psi_\tau$. After
recollecting some previous results in \cite{nek06,johnstone01} for
$\psi_\tau$, we give a detailed analysis of $\psi_\tau'$,
$\psi_\tau'-G'$ and also strengthen a previous bound on $\psi_\tau-
G$. Finally, we transfer the bounds on quantities related to $\psi
_\tau$
to those related to $\phi_\tau$ by a change of variable argument.

\subsection{Liouville--Green approach}
\label{subsec:louville-green}

Recall ($\tmu_{n,N},\tsigma_{n,N}$) in \eqref{eq:lag-center-scale} and
$\aa$ in \eqref{eq:ch-var}. We introduce the intermediate function
%
%
\begin{equation}\label{eq:fnN-def}
F_{n,N}(x) = (-1)^N \tsigma_{n,N}^{-1/2} \sqrt{N!/n!} x^{\aa/2 +
1}\mathrm{e}^{-x/2}L_N^{\aa+1}(x)
\end{equation}
as in \cite{johnstone01}, equation~(5.1), and \cite{nek06},  Section
2.2.2. (Note: $\aa=\aa_N - 1$ for the constant $\aa_N$
used in \cite{johnstone01} and \cite{nek06}.) Then $\phi_\tau$ is
related to $F_{n,N}$ as
\[
\psi_\tau(s) = \frac{1}{\sqrt{2}} \biggl(
\frac{N^{1/4}(n-1)^{1/4}\tsigma_{n-1,N-1}^{1/2}\sigma_{n,N}}{\tmu_{n-1,N-1}}
\biggr) F_{n-1,N-1}(\tauexpr{s}) \biggl(
\frac{\tmu_{n-1,N-1}}{\tauexpr{s}} \biggr).
\]
%
Replacing the subscripts $(n-1, N-1)$ by $(n-2, N)$ in $\tmu_{n-1,N-1},
\tsigma_{n-1,N-1}$ and $F_{n-1,N-1}$ on the right-hand side, we also obtain the
expression for $\phi_\tau(s)$.

Due to the close connection of $\psi_\tau$ and $\phi_\tau$ to $F_{n,N}$,
the key element in the proof of Proposition~\ref{prop:laguerre}
becomes asymptotic analysis of $F_{n,N}$ and its derivative. To this
end, the Liouville--Green (LG) theory set out in Olver \cite{olver74}, Chapter
11, is useful, for it comes with ready-made bounds on the
difference between $F_{n,N}$ and the Airy function, and also on the
difference between their derivatives.

To start with, we observe that $F_{n,N}$ satisfies a second-order
differential equation,
%
%
\begin{equation}\label{eq:diff-eqn}
F_{n,N}''(x)=\biggl \{\frac{1}{4}-\frac{\kk_N}{x}+\frac{\lambda_N^2-1/4}
{x^2} \biggr\} F_{n,N}(x),
\end{equation}
with $\kk_N = \frac{1}{2}(n+N+1)$ and $\lambda_N =
\frac{1}{2}(n-N)$. By rescaling $x = \kk_N\xi$, setting $w_N(\xi)
= F_{n,N}(x)$, the equation becomes
\[
w_N''(\xi) = \{\kappa_N^2 f(\xi)+g(\xi) \}w_N(\xi),
\]
where
\[
f(\xi) = \frac{(\xi- \xi_-)(\xi-\xi_+)}{4\xi^2},\qquad g(\xi) =
\frac{1}{4\xi^2}.
\]
The zeros of $f$ are given by $\xi_\pm= 2\pm\sqrt{4 - \ww_N^2}$ for
$\ww_N = 2\lambda_N/\kk_N$. They are called the turning points of the
differential equation, for each separates an interval in which the
solutions are oscillating from one in which they are of exponential
type. The LG approach introduces a new independent variable, $\zz$, and
dependent variable, $W$, as
\[
\zz\biggl(\frac{\mathrm{d}\zz}{\mathrm{d}\xi} \biggr)^2 = f(\xi),\qquad
W = \biggl(\frac{\mathrm{d}\zz}{\mathrm{d}\xi} \biggr)^{1/2}w_N.
\]
Then the differential equation takes the form $W''(\zz) =
\{\kk_N^2\zz+v(\ww_N, \zz) \}W(\zz)$. Without the perturbation
term $v(\ww_N,\zz)$, this is the Airy equation having linearly
independent solutions in terms of Airy functions $\Ai(\kk_N^{2/3}\zz)$
and $\mathrm{Bi}(\kk_N^{2/3}\zz)$. We focus on approximating the
recessive solution $\Ai(\kk_N^{2/3}\zz)$.

Let $\hat{f} = f/\zz$. \cite{olver74}, Theorem 11.3.1, gives that
\[
w_N(\xi) \propto\hat{f}^{-1/4}(\xi)\{\Ai(\kk_N^{2/3}\zz) +
\varepsilon_2(\kk_N, \xi)\},
\]
where, uniformly for $\xi\in[2,\infty)$, the error term $\varepsilon_2$
satisfies
%
\begin{eqnarray}\label
{eq:ee2-bd}
| \varepsilon_2(\kk_N, \xi) | & \leq&(\MM/\EE)(\kk
_N^{2/3}\zz)
\biggl[\exp\biggl\{\frac{\lambda_0}{\kk_N} F(\ww_N)\biggr\} - 1\biggr], \\\label{eq:ee2-d-bd}
| \partial_\xi\varepsilon_2(\kk_N, \xi) | & \leq&
\kk_N^{2/3}\hat{f}^{1/2}(\xi)(\NN/\EE)(\kk_N^{2/3}\zz)
\biggl[\exp\biggl\{\frac{\lambda_0}{\kk_N} F(\ww_N)\biggr\} - 1\biggr].
\end{eqnarray}
In the bounds, $\MM, \EE$ are the modulus and weight functions for the
Airy function and $\NN$ the phase function for its derivative (\cite
{olver74}, pages~394--396). On the real line, $\EE\geq1$ and
is increasing, $0\leq\MM\leq1$ and $\NN\geq0$. Moreover, for all $x$,
%
%
\begin{equation}
\label{eq:emn-ai}
|\Ai(x)|\leq(\MM/\EE)(x),\qquad|\Ai'(x)|\leq(\NN/\EE)(x).
\end{equation}
As $x\to\infty$, their asymptotics are given by
%
%
\begin{equation}
\label{eq:emn}
\EE(x)\sim\sqrt{2}\mathrm{e}^{ ({2}/{3})x^{3/2}}, \qquad\MM(x) \sim
\pi^{-1/2}x^{-1/4},\qquad\NN(x)\sim\pi^{-1/2}x^{1/4}.
\end{equation}
In addition, in the bounds \eqref{eq:ee2-bd} and \eqref{eq:ee2-d-bd},
$\lambda_0 \doteq1.04$ and the analysis in \cite{nek06}, A.3, shows
that, uniformly for $\xi\in[2,\infty)$, for large enough $N$,
%
%
\begin{equation}
\label{eq:F}
\exp\biggl\{\frac{\lambda_0}{\kk_N} F(\ww_N)\biggr\} - 1 \leq N^{-2/3}.
\end{equation}

Come back to $F_{n,N}$. The alignment in \cite{nek06}, equation (5)
and A.1,
shows that
\[
F_{n,N}(x) = r_N \kk_N^{1/6}\tsigma_{n,N}^{1/2}\hat{f}^{-1/4}(\xi)
\{\Ai(\kk_N^{2/3}\zz) + \varepsilon_2(\kk_N,\xi)\},
\]
with $r_N = 1 + \Oh{N^{-1}}$. Let $R_N(\xi) = (\zz'(\xi)/\zz'_N)^{-1/2}$
with $\zz'_N = \zz'(\xi_+)$. As $(\zz_N')^{-1} =
\kk_N^{1/3}\tsigma_{n,N}$ and $\hat{f}(\xi) = \zz'(\xi)^2$, we can
rewrite $F_{n,N}$ as
%
%
\begin{equation}
\label{eq:fnN}
F_{n,N}(x) = r_NR_N(\xi)\{\Ai(\kk_N^{2/3}\zz) +\varepsilon_2(\kk
_N, \xi)\}.
\end{equation}
This representation serves as the starting point for all the subsequent
asymptotic analysis on $\phi_\tau$, $\psi_\tau$ and their derivatives.

From now on, without notice, all the inequalities are understood to
hold uniformly for $N\geq N_0(s_0,\gamma)$.

\subsection{\texorpdfstring{Summary of previous analysis: Bound for $|\psi_\tau(s)|$}
{Summary of previous analysis: Bound for |psi tau(s)|}}
\label{subsec:psi}

Here, we summarize the previous analysis of $F_{n,N}$ in
\cite{johnstone01,nek06}, which gives the desired bound for
$|\psi_\tau(s)|$ in \eqref{eq:psi-psi'-bd} and a crude estimate for
$|\psi_\tau- G|$.

Let $x_{n,N}(s) = \xnexpr{s}$ and define
%
%
\begin{equation}
\label{eq:theta}
\theta_{n,N}(x_{n,N}(s)) =
F_{n,N}(x_{n,N}(s)) \biggl(\frac{\tmu_{n,N}}{x_{n,N}(s)} \biggr).
\end{equation}
As $\tsigma_{n,N}^{-1/2}N^{1/6} < 1$, we obtain that, for all $s\geq0$,
\[
|F_{n,N}(x_{n,N}(s))| \leq
|F_{n,N}(x_{n,N}(s))\tsigma_{n,N}^{1/2}N^{-1/6}|
\leq C\exp(-s),
\]
where the latter inequality was obtained in \cite{johnstone01}, A.8. If
$s_0<0$, then $\xi= x_{n,N}(s)/\kk_N\geq2$ uniformly for all $s\geq
s_0$.
In addition, Lemma \ref{lemma:I1N} later
shows that $|R_N(\xi)| \leq1 + CN^{-2/3}|s|$ for
$s\in[s_0,0]$. Therefore, we apply \eqref{eq:ee2-bd}, \eqref{eq:F} and
\eqref{eq:fnN} to obtain that
\[
| F_{n,N}(x_{n,N}(s)) | \leq2 r_N |R_N(\xi)|(\MM/\EE)(\kk
_N^{2/3}\zz)
\leq4,
\]
uniformly for $s\in[s_0,0]$.
Hence,
$|F_{n,N}(x_{n,N}(s))|\leq C \exp(-s)$ for all $s\geq s_0$. Moreover, we
note that $\tsigma_{n,N}/\tmu_{n,N} = \Oh{N^{-2/3}}$. So, when $N
\geq
N_0(s_0)$, for all $s\geq s_0$,
\[
{\tmu_{n,N}}/{x_{n,N}(s)} \leq(1 +
s_0 {\tsigma_{n,N}}/{\tmu_{n,N}} )^{-1} \leq2.
\]
Hence, uniformly for $s\geq s_0$, 
%
%
\begin{equation}
\label{eq:theta-bd}
| \theta_{n,N}(x_{n,N}(s)) | \leq C(s_0)\exp(-s).
\end{equation}
Finally, for any $\varrho_N = 1 + \Oh{N^{-1}}$, El Karoui \cite
{nek06}, Section
3.2, showed that, for all $s\geq s_0$,
\[
|\varrho_N\theta_{n,N}(x_{n,N}(s)) - \Ai(s)|
\leq C(s_0)N^{-2/3}\exp(-s/2).
\]

For $\psi_\tau(s)$, observe that $(\mu_{n,N},\sigma_{n,N}) =
(\tmu_{n-1,N-1}, \tsigma_{n-1,N-1})$. Using Sterling's formula, we
obtain that $\psi_\tau(s) = \frac{1}{\sqrt{2}} \rho_{N}
\theta_{n-1,N-1}(x_{n-1,N-1}(s))$ for some $\rho_N = 1 +
\Oh{N^{-1}}$. Then, we apply the last two displays to obtain
%
%
\begin{equation}
\label{eq:psi-tau-G-crude}
|\psi_\tau(s)|\leq C(s_0)\exp(-s),\qquad
| \psi_\tau(s) - G(s) | \leq C(s_0)N^{-2/3}\exp(-s/2),
\end{equation}
uniformly for $s\geq s_0$. 

Here, the first inequality gives the bound for $|\psi_\tau|$, while the
bound on $|\psi_\tau(s) - G(s)|$ could be further improved; see
\eqref{eq:psitau-g-improved}. Note that we cannot apply these results
directly to~$\phi_\tau$ since the `optimal' rescaling constants
$(\tmu_{n-2,N},\tsigma_{n-2,N})$ for $F_{n-2,N}$ do not agree with the
global constants $(\mu_{n,N}, \sigma_{n,N})$.

\subsection{\texorpdfstring{Asymptotics of $|\psi'_\tau(s)|$, $|\psi'_\tau(s)-G'(s)|$ and $|\psi_\tau(s)-G(s)|$}
{Asymptotics of |psi' tau(s)|, |psi' tau(s) - G'(s)| and |psi tau(s) - G(s)|}}
\label{sec:bd-psi}

Here, we derive bounds on $|\psi_\tau'|$ and $|\psi_\tau'-G'|$ and
refine the bound on $|\psi_\tau(s)-G(s)|$.

\subsubsection{\texorpdfstring{Bound for $|\psi'_\tau(s)|$}{Bound for |psi' tau(s)|}}

To obtain bounds for $|\psi_\tau'|$, we study
$|\partial_s\theta_{n,N}(x_{n,N}(s))|$. By the triangle inequality,
%
%
\begin{eqnarray}
\label{eq:d-theta-decomp}
| \partial_s \theta_{n,N}(x_{n,N}(s)) | & \leq&\biggl|
\tsigma_{n,N}F'_{n,N}(x_{n,N}(s))\frac{\tmu_{n,N}}{x_{n,N}(s)} \biggr|
+ \biggl|
\tsigma_{n,N}F_{n,N}(x_{n,N}(s))\frac{\tmu_{n,N}}{x^2_{n,N}(s)}
\biggr| \nonumber\qquad\quad
\\[-8pt]
\\[-8pt]
& =& T_{N,1}(s) + T_{N,2}(s).\qquad\quad
\nonumber
\end{eqnarray}
In what follows, we deal with the two terms in order.


\paragraph*{The $T_{N,1}$ term} Recall that $\tmu
_{n,N}/x_{n,N}(s)\leq
2$ for large $N$. So, we focus on $\tsigma_{n,N}F'_{n,N}$, which can be
decomposed as $\tsigma_{n,N}F'_{n,N} = \sum_{i=1}^4 D^i_{n,N}$, with
\begin{eqnarray*}
D^1_{n,N} & =& r_N\tsigma_{n,N}\kk_N^{-1}R'_N(\xi)
\{\Ai(\kk_N^{2/3}\xi) + \varepsilon_2(\kk_N, \xi)\}, \qquad
  D^2_{n,N}   =r_N[R_N^{-1}(\xi) - 1]\Ai'(\kk_N^{2/3}\zz),\\
D^3_{n,N} & =& r_N \Ai'(\kk_N^{2/3}\zz), \qquad
  D^4_{n,N} = r_N\tsigma_{n,N}\kk_N^{-1}R_N(\xi)
\partial_\xi\varepsilon_2(\kk_N, \xi).
\end{eqnarray*}
Due to different strategies used for the asymptotics on the $s$-scale,
we divide $[s_0,\infty)$ into $I_{1,N}\cup I_{2,N}$, with $I_{1,N} =
[s_0, s_1N^{1/6})$ and $I_{2,N} = [s_1N^{1/6}, \infty)$. The choice of
$s_1$ is worked out in Section~\ref{subsec:a-s-1}. Here, we note that $s_1
\geq1$ and that, for $s\geq s_1$,
%
%
\begin{equation}
\label{eq:EE-bd}
\EE^{-1}(\kk_N^{2/3}\zz) \leq C\exp(-3s/2) \leq C\exp(-s).
\end{equation}
In addition, we will repeatedly use the following facts.
\begin{lemma}
\label{lemma:I1N}
Under the conditions of Proposition \ref{prop:laguerre}, when $N\geq
N_0(s_0,\gamma)$, for all $s\in I_{1,N}$,
\begin{eqnarray*}
|R_N'(\xi)|&\leq& C\gamma^{-1/2}(1+\gamma),\qquad
|R_N(\xi)-1|\leq CN^{-2/3}|s|,\\
|\kk_N^{2/3}\zz- s| &\leq&(CN^{-2/3}s^2)\wedge\tfrac{1}{2}|s| \wedge1.
\end{eqnarray*}
\end{lemma}

Proof of Lemma \ref{lemma:I1N} is given in \cite{twacclong}.

\subparagraph*{Case $s\in I_{1,N}$} Consider $D^1_{n,N}$ first. Recall
that $r_N = 1 + \Oh{N^{-1}}$. Together with Lemma \ref{lemma:I1N}, this
implies
%
%
\begin{equation}
\label{eq:D-n-N-1-1}
|r_N\tsigma_{n,N}\kk_N^{-1}R_N'(\xi)| \leq CN^{-2/3}.
\end{equation}
On the other hand, as $0\leq\MM\leq1$, \eqref{eq:ee2-bd},
\eqref{eq:emn-ai} and \eqref{eq:F} together imply
\[
|\Ai(\kk_N^{2/3}\zz) + \varepsilon_2(\kk_N, \xi)| \leq
C (\MM/\EE)(\kk_N^{2/3}\zz)\leq C\EE^{-1}(\kk_N^{2/3}\zz).
\]
For $s\geq0$, Lemma \ref{lemma:I1N} implies $\kk_N^{2/3}\zz\geq
s/2$. Since $\EE$ is monotone increasing, by \eqref{eq:emn},
\[
|\Ai(\kk_N^{2/3}\zz) + \varepsilon_2(\kk_N, \xi)| \leq C\EE
^{-1}(s/2) \leq
C \mathrm{e}^{- ({1}/({3\sqrt{2}}))s^{3/2}}\leq C\exp(-s).
\]
If $s_0\leq0$, we can replace the $C$ on the rightmost side with
$C(s_0) = \max\{C,\break \max_{s\in[3s_0/2,0]}\EE^{-1}(s) \}$,
which is
continuous and non-increasing in $s_0$. Together with \eqref
{eq:D-n-N-1-1}, we
obtain that
\[
| D^1_{n,N} | \leq C(s_0)N^{-2/3}\exp(-s).
\]

(Here and after, we derive more stringent bounds
with the $N^{-2/3}$ term whenever possible. Although they are not
necessary for bounding $|\psi'_\tau|$, they are useful in the later study
of $|\psi'_\tau(s)-G'(s)|$.)

For $D_{n,N}^2$, we first have $|r_NR_N^{-1}(\xi) - 1| \leq
r_N|R_N^{-1}(\xi)-1| + |r_N - 1|$. Lemma \ref{lemma:I1N} implies that
$|R_N^{-1}(\xi)-1|\leq CN^{-2/3}|s|$. Observing that $|r_N-1| =
\Oh{N^{-1}}$, we obtain
\[
|r_NR_N^{-1}(\xi) - 1| \leq CN^{-2/3}|s|.
\]
For $|\Ai'(\kk_N^{2/3}\zz)|$, when $s\geq0$, Lemma \ref{lemma:I1N}
gives $\kk_N^{2/3}\zz\in[s/2, 3s/2]$. This, together with Lemma~\ref
{lemma:airy}, implies that
%
%
\begin{equation}
\label{eq:ai'-kk-bd}
|\Ai'(\kk_N^{2/3}\zz)|\leq C\exp(-3s/2).
\end{equation}
If $s_0 < 0$, we can replace the $C$ on the right-hand side with $C(s_0) =
\max\{C,\break \max_{[3s_0/2, 0]}|\Ai'(s)| \}$, which is continuous and
non-increasing.
Then the last two displays give
\[
|D^2_{n,N}| \leq C(s_0)N^{-2/3}|s|\exp(-3s/2)\leq
C(s_0)N^{-2/3}\exp(-s).\vspace*{-2pt}
\]

For $D^3_{n,N}$, we recall that $r_N = 1 + \Oh{N^{-1}}$. Together with
\eqref{eq:ai'-kk-bd}, this implies that
\[
|D^3_{n,N}|\leq C(s_0)\exp(-s).\vspace*{-2pt}
\]

For $D^4_{n,N}$, since $r_N = 1 +\Oh{N^{-1}}$, $\zz'(\xi) = \hat
{f}^{1/2}(\xi)$ and $\zz'_N = \kk_N^{1/3}/\tsigma_{n,N}$, \eqref
{eq:ee2-d-bd} and \eqref{eq:F} imply
\begin{eqnarray*}
|D^4_{n,N}| & =& |r_N \tsigma_{n,N}\kk_N^{-1} R_N(\xi)\partial_\xi
\varepsilon_2(\kk_N, \xi)| \\[-2pt]
& \leq& CN^{-2/3} \tsigma_{n,N} \kk_N^{-1/3} R_N(\xi) (\NN/\EE)(\kk
_N^{2/3}\zz)\\[-2pt]
& =& CN^{-2/3} R_N^{-1}(\xi)(\NN/\EE)(\kk_N^{2/3}\zz).\vspace*{-2pt}
\end{eqnarray*}
Lemma \ref{lemma:I1N} implies that $R_N^{-1}(\xi) \leq C$ and $\kk
_N^{2/3}\zz\in[s/2, 3s/2]$, uniformly on $I_{1,N}$. So, \eqref
{eq:emn} gives
\[
(\NN/\EE)(\kk_N^{2/3}\zz)\leq Cs^{1/4}\mathrm{e}^{- ({1}/({3\sqrt
{2}}))s^{3/2}} \leq C\exp(-s)\vspace*{-2pt}
\]
for all $s\geq0$. And if $s_0 < 0$, we can replace the $C$ on the
rightmost side with $C(s_0) = \max\{C, \max_{s\in[3s_0/2, 0]}(\NN
/\EE)(s)\}$, which is continuous and non-increasing in $s_0$. All
these elements together lead to
\[
|D^4_{n,N}| \leq C(s_0)N^{-2/3}\exp(-s).\vspace*{-2pt}
\]

Combining all the bounds on the $D^i_{n,N}$ terms, we obtain that $T_{N,1}
\leq C(s_0)\exp(-s)$ on $I_{1,N}$.

\subparagraph*{Case $s\in I_{2,N}$} In this case, we define
$\tilde{D}^1_{n,N} = D^1_{n,N}$ and $\tilde{D}^2_{n,N} = D^2_{n,N} +
D^3_{n,N} + D^4_{n,N}$.

Consider $\tilde{D}^1_{n,N}$ first. By \eqref{eq:ee2-bd},
\eqref{eq:emn-ai} and \eqref{eq:F}, we obtain that for $N\geq
N_0(s_0,\gamma)$,
\[
|\tilde{D}_{n,N}^1|\leq C\tsigma_{n,N}\kk_N^{-1}|R_N'/R_N|(\xi)
R_N(\xi)(\MM/\EE)(\kk_N^{2/3}\zz).\vspace*{-2pt}
\]
Observe that, uniformly on $I_{2,N}$,
%
%
\begin{equation}
\label{eq:tD-nN-1-bds}
\tsigma_{n,N}\kk_N^{-1}|R_N'/R_N|(\xi)\leq C,\qquad
R_N(\xi)\MM(\kk_N^{2/3}\zz) \leq Cs.\vspace*{-2pt}
\end{equation}
For a proof of \eqref{eq:tD-nN-1-bds}, see \cite{twacclong}.
On the other hand, \eqref{eq:EE-bd} holds on
$I_{2,N}$. Thus,
\[
|\tilde{D}_{n,N}^1| \leq Cs\exp(-3s/2) \leq Cs^4\exp(-s)
\leq CN^{-2/3}\exp(-s).\vspace*{-2pt}
\]

For $\tilde{D}^2_{n,N}$, we can write it as $\tilde{D}^2_{n,N} = r_N
R_N(\xi)[\Ai'(\kk_N^{2/3}\zz)R_N^{-2}(\xi) + \tsigma_{n,N}\kk_N^{-1}
\partial_\xi\varepsilon_2(\kk_N, \xi)]$. By \eqref{eq:ee2-d-bd},
\eqref{eq:emn-ai} and \eqref{eq:F} and the identity
$R_N^{-1}=\tsigma_{n,N}^{-1/2}\kk_N^{1/6}\hat{f}^{1/4}$, we get the
bound
\[
|\tilde{D}_{n,N}^2| \leq C
R_N^{-1}(\xi)(\NN/\EE)(\kk_N^{2/3}\zz).\vspace*{-2pt}\vadjust{\goodbreak}
\]
\eqref{eq:emn} suggests that $R_N^{-1}(\xi)\NN(\kk_N^{2/3}\zz)\leq
CR_N^{-1}(\xi)\kk_N^{1/6}\zz^{1/4} =
Cf^{1/4}(\xi)\tsigma_{n,N}^{1/2}\leq C\tsigma_{n,N}^{1/2}$. The last
inequality holds as $f\leq4$ for $s\in I_{2,N}$. On the other hand,
$\tsigma_{n,N}\leq C(\gamma)N^{1/3}\leq Cs^4$ for large $N$. Assembling
all the pieces, we obtain $R_N^{-1}(\xi)\NN(\kk_N^{2/3}\zeta) \leq
Cs^2$. Together with~\eqref{eq:EE-bd}, this implies
\[
|\tilde{D}^2_{n,N}| \leq Cs^2\exp(-3s/2) \leq
Cs^{-4}\exp(-s)\leq CN^{-2/3}\exp(-s).
\]
Therefore, $T_{N,1}\leq CN^{-2/3}\exp(-s)$ on $I_{2,N}$.

\paragraph*{The $T_{N,2}$ term}
This term is relatively easy to bound. Note that
$\tsigma_{n,N}/\tmu_{n,N} = \Oh{N^{-2/3}}$ and that $T_{N,2}(s) =
|\theta_{n,N}(x_{n,N}(s))\tsigma_{n,N}/x_{n,N}(s)|$. So, for all
$s\geq s_0$,
$N\geq N_0(s_0)$,
\[
|\tsigma_{n,N}/x_{n,N}(s)| = |s + \tmu_{n,N}/\tsigma_{n,N}|^{-1}
\leq
C(s_0)N^{-2/3}.
\]
Together with \eqref{eq:theta-bd}, this implies that for all $s\geq
s_0$, $T_{N,2}(s)\leq C(s_0)N^{-2/3}\exp(-s)$.

\paragraph*{Summing up} By \eqref{eq:d-theta-decomp}, the bounds on
$T_{N,1}$ and $T_{N,2}$ transfer to
%
%
\begin{equation}
\label{eq:d-theta}
| \partial_s\theta_{n,N}(x_{n,N}(s)) | \leq
C(s_0)\exp(-s)
\end{equation}
uniformly for $s\geq s_0$. On the other
hand, we note that
\[
\psi'_\tau(s) = \tfrac{1}{\sqrt{2}}\rho_N
\partial_s\theta_{n-1,N-1}(x_{n-1,N-1}(s)),
\]
with $\rho_N = 1 + \Oh{N^{-1}}$. Thus, \eqref{eq:d-theta} implies the
desired bound on $|\psi'_\tau|$ in \eqref{eq:psi-psi'-bd}.

\subsubsection{\texorpdfstring{Bound for $|\psi'_\tau(s)-G'(s)|$}{Bound for |psi' tau(s) - G'(s)|}}

By the triangle inequality, we bound $|\psi'_\tau(s)-G'(s)|$ as
%
%
\begin{eqnarray}
\label{eq:psi'-tau-G'}
|\psi'_\tau(s)-G'(s)| &\leq& \tfrac{1}{\sqrt{2}}|\rho_N-1|
|\partial_s\theta_{n-1,N-1}(x_{n-1,N-1}(s))| \nonumber
\\[-8pt]
\\[-8pt]
&&{}
+ \tfrac{1}{\sqrt{2}}
|\partial_s\theta_{n-1,N-1}(x_{n-1,N-1}(s))-\Ai'(s)|.
\nonumber
\end{eqnarray}
As $\rho_N = 1 + \Oh{N^{-1}}$, by \eqref{eq:d-theta}, we bound the first
term by $C(s_0)N^{-1}\exp(-s)$. In what follows, to bound the second
term in \eqref{eq:psi'-tau-G'}, we focus on
$|\partial_s\theta_{n,N}(x_{n,N}(s))-\Ai'(s)|$, which can first be split
into two parts as:
\begin{eqnarray*}
&& |\partial_s\theta_{n,N}(x_{n,N}(s))-\Ai'(s)| \\
&& \quad  \leq\biggl|
\tsigma_{n,N}F'_{n,N}(x_{n,N}(s))\frac{\tmu_{n,N}}{x_{n,N}(s)} -
\Ai'(s)\biggr| + \biggl|
\tsigma_{n,N}F_{n,N}(x_{n,N}(s))\frac{\tmu_{n,N}}{x_{n,N}^2(s)}
\biggr|
\\
&& \quad  = \mathcal{T}_{N,1}(s) + \mathcal{T}_{N,2}(s).
\end{eqnarray*}

\paragraph*{The $\mathcal{T}_{N,1}(s)$ term}
For this term, we separate the arguments on\vspace*{1pt} $I_{1,N}=[s_0,s_1N^{1/6})$
and $I_{2,N} = [s_1N^{1/6},\infty)$.

\subparagraph*{Case $s\in I_{1,N}$} On $I_{1,N}$, we decompose
$\mathcal{T}_{N,1}(s)$ as $ \mathcal{T}_{N,1}(s) =
\sum_{i=1}^5\mathcal{D}^i_{n,N}$, with $\mathcal{D}^i_{n,N} =
D^i_{n,N}\tmu_{n,N}/x_{n,N}(s)$ for $i = 1, 2$ and $4$, and
\[
\mathcal{D}^3_{n,N} = r_N\frac{\tmu_{n,N}}{x_{n,N}(s)}
[\Ai'(\kk_N^{2/3}\zz) - \Ai'(s)],\qquad
\mathcal{D}^5_{n,N} = \biggl[r_N\frac{\tmu_{n,N}}{x_{n,N}(s)}-1 \biggr]
\Ai'(s).\vspace*{-2pt}
\]
Observe that $|\tmu_{n,N}/x_{n,N}(s)|\leq2$ on $I_{1,N}$. Thus, by
previous bounds on $D_{n,N}^i$, we obtain that, for $i=1,2$ and $4$,
$|\mathcal{D}^i_{n,N}| \leq C(s_0)N^{-2/3}\exp(-s)$.

Consider $\mathcal{D}^3_{n,N}$. By the Taylor expansion, for some $s^*$
between $\kk_N^{2/3}\zz$ and $s$,
\[
|\Ai'(\kk_N^{2/3}\zz) - \Ai'(s)| \leq
|\Ai''(s^*)| |\kk_N^{2/3}\zz- s| =
|s^*\Ai(s^*)| |\kk_N^{2/3}\zz- s|,\vspace*{-2pt}
\]
where the equality comes from the identity $\Ai''(s) = s\Ai(s)$. By
Lemma \ref{lemma:I1N}, we have that $|\kk_N^{2/3}\zz- s| \leq CN^{-2/3}s^2$,
and that $s^*$ lies between $\frac{1}{2}s$ and $\frac{3}{2}s$. The
latter, together with Lemma \ref{lemma:airy}, implies that, for $s\geq0$,
\[
|s^*\Ai(s^*)| \leq C\exp(-3s/2).\vspace*{-2pt}
\]
If $s_0\leq0$, we then have $s^*\in[\frac{3}{2}s, 0]$, and hence we
can replace $C$ on the right-hand side with $C(s_0) = \max\{C,\max_{s\in
[3s_0/2,0]} |s\Ai(s)|\}$. Observe that $r_N = 1 + \Oh{N^{-1}}$ and
that $|\tmu_{n,N}/x_{n,N}(s)|\leq2$. We thus conclude that
\[
|\mathcal{D}^3_{n,N}| \leq C(s_0)N^{-2/3}s^2\exp(-3s/2)
\leq C(s_0)N^{-2/3}\exp(-s).\vspace*{-2pt}
\]

Switch to $\mathcal{D}^5_{n,N}$. We first note that
\begin{eqnarray*}
\biggl| r_N\frac{\tmu_{n,N}}{x_{n,N}(s)}-1\biggr| &\leq& r_N\biggl|
\frac{\tmu_{n,N}}{x_{n,N}(s)}-1\biggr| + | r_N - 1 | \\[-2pt]
& =& r_N|s|\biggl| s+ \frac{\tmu_{n,N}}{\tsigma_{n,N}} \biggr| ^{-1} +
| r_N-1 | \leq CN^{-2/3}|s|+CN^{-1}.\vspace*{-2pt}
\end{eqnarray*}
The last inequality holds as $\tsigma_{n,N}/\tmu_{n,N}=\Oh
{N^{-2/3}}$, $r_N = 1 + \Oh{N^{-1}}$, and for large~$N$, $|s +
\tmu_{n,N}/\tsigma_{n,N}|\geq\frac{1}{2}\tmu_{n,N}/\tsigma_{n,N}$
uniformly for $s\in I_{1,N}$. On
the other hand, Lemma \ref{lemma:airy} implies that $|\Ai'(s)|\leq
C(s_0)\exp(-3s/2 )$.
Putting the two parts together, we obtain
\[
|\mathcal{D}^5_{n,N}| \leq C(s_0)N^{-2/3}(|s| + CN^{-1/3})
\exp(-3s/2) \leq C(s_0)N^{-2/3}\exp(-s).\vspace*{-2pt}
\]

Assembling all the bounds on the $\mathcal{D}^i_{n,N}$'s, we obtain
that, on $I_{1,N}$,
\[
\mathcal{T}_{N,1}(s)\leq C(s_0)N^{-2/3}\exp(-s).\vspace*{-2pt}
\]

\subparagraph*{Case $s\in I_{2,N}$} In this case, we could act more
heavy-handedly. In particular, by the asymptotics of $T_{N,1}(s)$ on
$I_{2,N}$ and Lemma \ref{lemma:airy}, we have
\begin{eqnarray*}
\mathcal{T}_{N,1}(s) & \leq&\biggl|
\tsigma_{n,N}F'_{n,N}(x_{n,N}(s))\frac{\tmu_{n,N}}{x_{n,N}(s)} \biggr|
+ | \Ai'(s)|
\leq CN^{-2/3}\exp(-s) + C\exp(-3s/2) \\
& \leq& CN^{-2/3}\exp(-s) + CN^{-2/3}s^{4}\exp(-3s/2) \leq
CN^{-2/3}\exp(-s).\vspace*{-2pt}\vadjust{\goodbreak}
\end{eqnarray*}

\paragraph*{The $\mathcal{T}_{N,2}(s)$ term}
The $\mathcal{T}_{N,2}(s)$ term is the same as $T_{N,2}(s)$ defined
previously in the study of $\partial_s\theta_{n,N}(x_{n,N}(s))$ and
hence we quote the bound derived there directly as
\[
\mathcal{T}_{N,2}(s) \leq C(s_0)N^{-2/3}\exp(-s) \qquad
\mbox{for all }s\geq s_0.
\]

\paragraph*{Summing up} Combining the bounds on $\mathcal{T}_{N,1}$ and
$\mathcal{T}_{N,2}$, we have, uniformly for $s\geq s_0$,
\[
|\partial_s\theta_{n,N}(x_{n,N}(s))-\Ai'(s)| \leq
C(s_0)N^{-2/3}\exp(-s).
\]
By the discussion following \eqref{eq:psi'-tau-G'}, we obtain the
desired bound on $|\psi'_\tau(s)-G'(s)|$ in \eqref{eq:psi'-bd}.

\subsubsection{\texorpdfstring{Improved bound for $|\psi_\tau-G|$}{Improved bound for |psi tau - G|}}
The bound on $|\psi'_\tau(s)-G'(s)|$, together with
\eqref{eq:psi-tau-G-crude}, can lead to a tighter bound for
$|\psi_\tau(s)-G(s)|$ as the following:
%
%
\begin{eqnarray}\label{eq:psitau-g-improved}
| \psi_\tau(s) -G(s)| & =& \biggl|
\int_s^{2s}[\psi'_\tau(t)-G'(t)]\,\mathrm{d}t - [ \psi_\tau(2s) - G(2s)
] \biggr| \nonumber\\
& \leq&\int_s^{2s}| \psi'_\tau(t) - G'(t) |\,\mathrm{d}t + |
\psi_\tau(2s)-G(2s) | \\
& \leq&\int_s^{2s}C(s_0)N^{-2/3}\mathrm{e}^{-t}\,\mathrm{d}t +
C(s_0)N^{-2/3}\exp(-s)
\leq C(s_0)N^{-2/3}\exp(-s).\nonumber
\end{eqnarray}
This is exactly what we claimed in Proposition \ref{prop:laguerre}.

\subsection{\texorpdfstring{Asymptotics for quantities related to $\phi_\tau(s)$}
{Asymptotics for quantities related to phi tau(s)}}
\label{subsec:phi-bd}

In this part, we employ a trick in \cite{johnstone01} to transfer the
bounds on the quantities related to~$\psi_\tau$ to those related to
$\phi_\tau$.

Recall that, for $\tilde{\rho}_N = 1+\Oh{N^{-1}}$ (see Section \ref
{subsec:a-constants} for its proof),
\[
\phi_\tau(s) =
\frac{1}{\sqrt{2}}\tilde{\rho}_NF_{n-2,N}(x_{n-1,N-1}(s))\frac
{\tmu_{n-2,N}}
{x_{n-1,N-1}(s)}.
\]
If the $x_{n-1,N-1}(s)$ term on the right-hand side were $x_{n-2,N}(s)$,
then all the bounds we have proved for $\psi_\tau$ would also be
valid for
$\phi_\tau$. As this is not the case, we introduce a~new independent
variable $s'$ as{}:
%
%
\begin{equation}
\label{eq:s-prime-def}
x_{n-1,N-1}(s) = x_{n-2,N}(s'),
\end{equation}
that is, $s'=(\tmu_{n-1,N-1}-\tmu_{n-2,N})/\tsigma_{n-2,N} +
s\tsigma_{n-1,N-1}/\tsigma_{n-2,N}$. (The readers are expected not to
confuse it
with the $s' $ that previously appeared in Section \ref{sec:det-law}.)
Then $\phi_\tau$ can be rewritten as
\[
\phi_\tau(s)=\frac{1}{\sqrt{2}}\tilde{\rho}_NF_{n-2,N}(x_{n-2,N}(s'))
\frac{\tmu_{n-2,N}}{x_{n-2,N}(s')} = \frac{1}{\sqrt{2}}\tilde{\rho}_N
\theta_{n-2,N}(x_{n-2,N}(s')).
\]
Recalling the definition of $\Delta_N$ in \eqref{eq:Delta-N}, we have
$s'-s = \Delta_N + [\tsigma_{n-1,N-1}\tsigma_{n-2,N}^{-1}]s$, with
%
%
\begin{equation}
\label{eq:Delta-N-bd}
\Delta_N = \Oh{N^{-1/3}},\qquad
1\leq{\tsigma_{n-1,N-1}}{\tsigma_{n-2,N}^{-1}}=1+\Oh{N^{-1}}.
\end{equation}

\subsubsection*{Bounds for $|\phi_\tau(s)|$ and $|\phi'_\tau(s)|$}
Recall previous bounds on $|\theta_{n,N}(x_{n,N}(s))|$ and
$|\partial_s\theta_{n,N}(x_{n,N}(s))|$. Together with
\eqref{eq:Delta-N-bd}, they imply that, for all $s\geq s_0$,
\[
| \phi_\tau(s)| \leq C(s_0)\exp(-s')\leq C(s_0)\exp(-s)
\]
and
\begin{eqnarray*}
| \phi'_\tau(s)| &=& \frac{1}{\sqrt{2}}\tilde{\rho}_N|
\partial_s\theta_{n-2,N}(x_{n-2,N}(s'))|\\
 &=&
\frac{1}{\sqrt{2}}\tilde{\rho}_N| \partial_{s'}
\theta_{n-2,N}(x_{n-2,N}(s'))| \frac{\mathrm{d}s'}{\mathrm{d}s} \\
& \leq& C(s_0)\exp(-s')\frac{\tsigma_{n-1,N-1}}{\tsigma_{n-2,N}}
\leq C(s_0)\exp(-s).
\end{eqnarray*}

\subsubsection*{Bounds for $|\phi_\tau(s)-G_N(s)|$ and
$|\phi'_\tau(s)-G_N'(s)|$}

We consider $|\phi_\tau(s)-G_N(s)|$ in detail and the derivation for the
bound on $|\phi'_\tau(s)-G_N'(s)|$ is essentially the same.

By the definition of $s'$ and the identity $\Ai''(s) = s\Ai(s)$, we
obtain the Taylor expansion
\begin{eqnarray*}
G(s') & =& G(s) + (s'-s)G'(s) + \frac{1}{2}(s'-s)^2G''(s^*) \\
& =& G_N(s) +
\frac{1}{\sqrt{2}} \biggl[\frac{\tsigma_{n-1,N-1}}{\tsigma_{n-2,N}}-1
\biggr]s\Ai'(s) + \frac{1}{2\sqrt{2}}(s'-s)^2s^*\Ai(s^*),
\end{eqnarray*}
with $s^*$ lying in between $s$ and $s'$. By the previous discussion on
$|\psi_\tau(s)-G(s)|$, this leads to
%
%
\begin{eqnarray}
\label{eq:phi-tau-G}
|\phi_\tau(s) - G_N(s)| & \leq& C(s_0)N^{-2/3}\exp(-s')
+CN^{-1}|s\Ai'(s)| + C(s'-s)^2| s^*\Ai(s^*)|\nonumber\qquad
\\[-8pt]
\\[-8pt]
& \leq& C(s_0)N^{-2/3}\exp(-s) + C(s'-s)^2| s^*\Ai(s^*)|.
\nonumber\qquad
\end{eqnarray}
To further bound the last term, we split $[s_0,\infty)$ into
$I_{1,N}\cup I_{2,N}$. For $s\in I_{1,N}$,
\[
(s-s')^2=\biggl[\Delta_N + \biggl(\frac{\tsigma_{n-1,N-1}}{\tsigma_{n-2,N}}-1\biggr)s\biggr]^2
\leq[CN^{-1/3} + CN^{-1}s]^2 \leq(CN^{-2/3})\wedge1.
\]
So $|s^*|\leq|s|+ 1$, and Lemma \ref{lemma:airy} implies that
\[
C(s-s')^2|s^*\Ai(s^*)| \leq C(s_0)N^{-2/3}\exp(-s).
\]
On $I_{2,N}$, \eqref{eq:Delta-N-bd} implies that $s'\geq s/2$, and hence
$s^*\geq s/2$. Together with Lemma \ref{lemma:airy}, this implies
\[
C(s'-s)^2 |s^*\Ai(s^*)| \leq Cs^{-4}\cdot|(s^*)^7\Ai(s^*)|
\leq CN^{-2/3}\exp(-s).
\]
Therefore, we have shown that, for all $s\geq s_0$, the last term in
\eqref{eq:phi-tau-G} is further controlled by $C(s_0)N^{-2/3}\exp(-s)$,
which in turn gives the desired bound for $|\phi_\tau-G_N|$. It is not
hard to check that all the $C(s_0)$ functions in the above analysis
could be continuous and non-increasing.



\begin{appendix}\label{appm}
\section*{Appendix}
\label{appendix}
 \renewcommand{\theequation}{\arabic{equation}}

In the Appendix, we collect technical details that led to some of the
claims previously made in the main text. Section \ref{subsec:a-constants}
gives proofs
to properties of a number of constants. Section \ref{subsec:a-s-1} works out
the details on the choice of $s_1$, which was used to decompose the
interval $[s_0,\infty)$ in Section~\ref{sec:lg}.

\subsection{\texorpdfstring{Properties of $\beta_N, \rho_N, \tilde{\rho}_N, \Delta_N$
and $\tsigma_{n-1,N-1}/\tsigma_{n-2,N}$}
{Properties of beta N, rho N, tilde{rho} N, Delta N
and tilde{sigma} {n-1,N-1}/tilde{sigma {n-2,N}}}}
\label{subsec:a-constants}

\subsubsection*{Property of $\beta_N$}
We are to show that $\beta_N = \frac{1}{\sqrt{2}} + \Oh{N^{-1}}$. By
definition, we know
\begin{eqnarray*}
\beta_N & =& \frac{1}{2}\int_{-\infty}^\infty\phi_\tau(s)\,\mathrm{d}s =
\frac{1}{2}\int_0^\infty\phi(x;\aa)\,\mathrm{d}x \\
& =&
\frac{N^{1/4}(n-1)^{1/4}\Gamma^{1/2}(N+1)}{2\sqrt{2}\Gamma
^{1/2}(n)}\times
\int_0^\infty x^{(\aa-1)/2}\mathrm{e}^{-x/2}L_N^{\aa}(x)\,\mathrm{d}x \\
& =& \frac{2^{-\aa/2}N^{1/4}(n-1)^{1/4}\Gamma^{1/2}(n)\Gamma(
({1}/{2})(N+3))}
{(N+1)\Gamma^{1/2}(N+1)\Gamma( ({1}/{2})(n+1))}.
\end{eqnarray*}
Applying Sterling's formula $ \Gamma(z) =
({2\uppi}/{z} )^{1/2} ({z}/{\mathrm{e}} )^{z}(1+\Oh{z^{-1}})$,
we obtain that
\begin{eqnarray*}
\beta_N & =& \frac{ ({2\uppi}/{n} )^{1/4}
({n}/{\mathrm{e}} )^{n/2}
[{4\uppi}/(N+3) ]^{1/2} [(N+3)/(2\mathrm{e})
]^{(N+3)/2}}{ [{2\uppi}/(N+1) ]^{1/4} [(N+1)/{\mathrm{e}}
]^{(N+1)/2} [{4\uppi}/(n+1) ]^{1/2} [(n+1)/(2\mathrm{e})
]^{(n+1)/2}} \\
&&{} \times\frac{2^{-\aa/2}N^{1/4}(n-1)^{1/4}}{N+1}\bigl (1 + \Oh
{N^{-1}}\bigr) \\
& =& \frac{1}{\sqrt{2\mathrm{e}}} \biggl(1 - \frac{1}{n+1} \biggr)^{n/2} \biggl(1
+ \frac{2}{N+1} \biggr)^{(N+1)/2 + 3/4} \bigl(1 + \Oh{N^{-1}}\bigr) \\
& =& \frac{1}{\sqrt{2}} + \Oh{N^{-1}}.
\end{eqnarray*}

\subsubsection*{Properties $\rho_N$ and $\tilde{\rho}_N$}
We want to show that $\rho_N, \tilde{\rho}_N = 1 + \Oh{N^{-1}}$.
Consider $\rho_N$ first. By definition, we have
\[
\rho_N =
\frac{N^{1/4}(n-1)^{1/4}\tsigma_{n-1,N-1}^{1/2}\sigma_{n,N}}{\mu_{n,N}}
= \frac{N^{1/4}(n-1)^{1/4}\tsigma_{n-1,N-1}^{3/2}}{\tmu_{n-1,N-1}}.
\]
Plugging in the definition of $\tsigma_{n-1,N-1}$ and $\tmu_{n-1,N-1}$,
we obtain that
\begin{eqnarray*}
\rho_N & =& N^{1/4}(n-1)^{1/4}
\biggl(\sqrt{N-\frac{1}{2}}+\sqrt{n-\frac{1}{2}} \biggr)^{-1/2}
\biggl(\frac{1}{\sqrt{N- {1}/{2}}}+\frac{1}{\sqrt{n-
{1}/{2}}} \biggr)^{1/2}
\\[-2pt]
& =& \biggl(\frac{N}{N- {1}/{2}} \biggr)^{1/4}
\biggl(\frac{n-1}{n- {1}/{2}} \biggr)^{1/4}
= 1 + \Oh{N^{-1}}.
\end{eqnarray*}

For $\tilde{\rho}_N$, we have
\begin{eqnarray*}
\tilde{\rho}_N & =& \frac{N^{1/4}(n-1)^{1/4}\tsigma_{n-2,N}^{1/2}
\sigma_{n,N}}{\tmu_{n-2,N}} =
\frac{\tsigma_{n-1,N-1}}{\tsigma_{n-2,N}}
\frac{N^{1/4}(n-1)^{1/4}\sigma_{n-2,N}^{3/2}}{\mu_{n-2,N}}\\[-2pt]
& =& \frac{\tsigma_{n-1,N-1}}{\tsigma_{n-2,N}}N^{1/4}(n-1)^{1/4}
\biggl(\sqrt{N+\frac{1}{2}}+\sqrt{n-\frac{3}{2}} \biggr)^{-1/2}
\biggl(\frac{1}{\sqrt{N+ {1}/{2}}}+\frac{1}{\sqrt{n-
{3}/{2}}} \biggr)^{1/2}\\[-2pt]
& =& \frac{\sigma_{n-1,N-1}}{\sigma_{n-2,N}}
\biggl(\frac{N}{N+ {1}/{2}} \biggr)^{1/4}
\biggl(\frac{n-1}{n- {3}/{2}} \biggr)^{1/4} = 1 + \Oh{N^{-1}}.
\end{eqnarray*}
The last equality holds since $\tsigma_{n-1,N-1}/\tsigma_{n-2,N} = 1 +
\Oh{N^{-1}}$ as claimed in \eqref{eq:Delta-N}, which is to be shown
below.

\subsubsection*{Property of $\Delta_N$} Recall the definition
$\Delta_N = (\tmu_{n-1,N-1}-\tmu_{n-2,N})/{\tsigma_{n-2,N}}$. By
\cite{nek06}, A.1.2, the numerator $\tmu_{n-1,N-1} -
\tmu_{n-2,N} = \Oh{1}$. For the denominator, let
$\gamma_{n,N} = (n-\frac{3}{2} )/ (N+\frac
{1}{2} )$. We
then have
\begin{eqnarray*}
\frac{1}{\tsigma_{n-2,N}} & =&
\biggl(\sqrt{N+\frac{1}{2}}+\sqrt{n-\frac{3}{2}} \biggr)^{-1}
\biggl(\frac{1}{\sqrt{N+ {1}/{2}}}+\frac{1}{\sqrt{n-
{3}/{2}}} \biggr)^{-1/3}
\\[-2pt]
& =& \frac{1}{1 + \gamma_{n,N}^{1/2}} (1 +
{\gamma_{n,N}^{-1/2}} ) \biggl(N +
\frac{1}{2} \biggr)^{-1/3} = \Oh{N^{-1/3}}.
\end{eqnarray*}
The last equality holds since $\gamma_{n,N}$ is bounded below for
all $n > N$. Combining the two parts, we establish that
$\Delta_{N} = \Oh{N^{-1/3}}$.

\subsubsection*{Property of $\tsigma_{n-1,N-1}/\tsigma_{n-2,N}$}
We now switch to prove that
\[
1 \leq{\tsigma_{n-1,N-1}}/{\tsigma_{n-2,N}} = 1 + \Oh{N^{-1}}.\vadjust{\goodbreak}
\]
\cite{nek06}, A.1.3, showed that $\tsigma_{n-1,N-1}/\tsigma_{n-2,N}
= 1 + \Oh{N^{-1}}$. On the other hand, we have
from the second-to-last display of \cite{nek06}, A.1.3, that
\[
\biggl( \frac{\tsigma_{n-1,N-1}}{\tsigma_{n-2,N}} \biggr)^3 =
\biggl[ 1
+ \frac{\sqrt{n/N} - \sqrt{N/n}}{n+N} + \Oh{n^{-2}} \biggr]\biggl [1
+ \frac{1}{2} \biggl(\frac{1}{n}+\frac{1}{N} \biggr) +
\Oh{n^{-2}} \biggr].
\]
Both terms become greater than $1$ when $N\geq N_0(\gamma)$, and hence
$\tsigma_{n-1,N-1}/\tsigma_{n-2,N} \geq1$ for large~$N$. Actually,
the inequality holds for any $n>N\geq2$. However, what we have proved
here is sufficient for our argument in Section \ref{subsec:phi-bd}.

\subsection{Choice of $s_1$ and its consequences}
\label{subsec:a-s-1}

The key point in our choice of $s_1$ is to ensure that when $s\geq
s_1$, we have
%
%
\begin{equation}
\label{eq:kk-N-zz-s-1-1}
\tfrac{2}{3}\kk_N\zz^{3/2}\geq\tfrac{3}{2}s.
\end{equation}

To this end, recall that in \cite{johnstone01}, A.8, one could choose
$\tilde{s}_1(\gamma) = C(\gamma)(1+\delta)$ with some $\delta>0$, such
that when $s\geq\tilde{s}_1(\gamma)$, we have $\sqrt{f(\xi)}\geq
2/\tsigma_{n,N}$ and hence if $s\geq4\tilde{s}_1(\gamma)$,
\[
\frac{2}{3}\kk_N\zz^{3/2} = \kk_N\int_{\xi_+}^\xi\sqrt{f(z)}\,\mathrm{d}z
\geq
\kk_N\frac{2}{\tsigma_{n,N}}\bigl(s-\tilde{s}_1(\gamma)\bigr)\frac{\tsigma
_{n,N}}{\kk_N}
= 2\bigl(s-\tilde{s}_1(\gamma)\bigr)\geq\frac{3}{2}s.
\]
Moreover, by the analysis in \cite{nek06}, A.6.4,
$\tilde{s}_1(\gamma)$ could be chosen independently of $\gamma$ and
hence we could define our $s_1$ to be
\[
s_1 = 4\tilde{s}_1,
\]
which is independent of $\gamma$ and such that \eqref{eq:kk-N-zz-s-1-1}
holds. Moreover, we also require that $s_1 \geq1$.

After specifying our choice of $s_1$, we spell out two of its
consequences. The first of them is that when $s\geq
s_1\geq1$,
%
%
\begin{equation}
\label{eq:EE-bd-1}
\EE^{-1}(\kk_N^{2/3}\zz) \leq C\exp(-3s/2) \leq C\exp(-s).
\end{equation}
This is from the observation that $\EE(x)\geq
C\exp(2 x^{3/2}/3)$ and hence
\[
\EE^{-1}(\kk_N^{2/3}\zz)\leq
C\exp\biggl(-\frac{2}{3}\kk_N\zz^{3/2} \biggr)\leq C\exp(-3s/2).
\]

The other consequence is about the behavior of $s'$ defined in
\eqref{eq:s-prime-def} when $s\geq s_1$. Remembering that $s_1\geq1$,
we then have that when $s\geq s_1$ and $N\geq N_0(\gamma)$,
%
%
\begin{equation}
\label{eq:s-prime-s/2}
s'-\frac{s}{2} = \Delta_N + \biggl(\frac{\tsigma_{n-1,N-1}}{\tsigma
_{n-2,N}} -
\frac{1}{2} \biggr)s\geq\Delta_N + \frac{s_1}{2} \geq\Delta_N +
\frac{1}{2} \geq0.
\end{equation}
The last inequality holds when $N\geq N_0(\gamma)$, for $\Delta_N =
\Oh{N^{-1/3}}$.
\end{appendix}

\section*{Acknowledgements}
I am most grateful to Professor Iain Johnstone for numerous
discussions. Thanks also go to Professor Debashis Paul for
kindly sharing an unpublished manuscript. I am grateful to the editor,
an associate editor and an anonymous referee for their helpful comments
that led to improvement on the presentation of the paper. This work is
supported in part by Grants NSF DMS-05-05303 and NIH EB R01 EB001988.

%

\printhistory

\end{document}